\documentclass[11pt]{article}
\usepackage{latexsym,amsfonts,amsmath,amssymb,amsthm,url}
\usepackage[english]{babel}
\usepackage{xcolor}

\usepackage{chngcntr}
\counterwithin{equation}{section}

\newtheorem{itlemma}{Lemma}[section]
\newtheorem{itproposition}[itlemma]{Proposition}
\newtheorem{theorem}[itlemma]{Theorem}
\newtheorem{itcorollary}[itlemma]{Corollary}
\newtheorem{itremark}[itlemma]{Remark}
\newtheorem{itremarks}[itlemma]{Remarks}
\newtheorem{itdefinition}[itlemma]{Definition}
\newtheorem{itexample}[itlemma]{Example}

\newenvironment{fact}{\begin{itfact}\rm}{\end{itfact}}

\newenvironment{lemma}{\begin{itlemma}}{\end{itlemma}}
\newenvironment{remark}{\begin{itremark}\rm}{\end{itremark}}
\newenvironment{remarks}{\begin{itremarks} \rm}{\end{itremarks}}
\newenvironment{corollary}{\begin{itcorollary}}{\end{itcorollary}}
\newenvironment{proposition}{\begin{itproposition}}{\end{itproposition}}
\newenvironment{definition}{\begin{itdefinition}\rm}{\end{itdefinition}}
\newenvironment{example}{\begin{itexample}\rm}{\end{itexample}}
\newcommand{\bl}[1]{\begin{lemma}\label{#1}}
\newcommand{\br}[1]{\begin{remark}\label{#1}}
\newcommand{\brs}[1]{\begin{remarks}\label{#1}}
\newcommand{\bt}[1]{\begin{theorem}\label{#1}}
\newcommand{\bd}[1]{\begin{definition}\label{#1}}
\newcommand{\bp}[1]{\begin{proposition}\label{#1}}
\newcommand{\bc}[1]{\begin{corollary}\label{#1}}
\newcommand{\bfact}[1]{\begin{fact}\label{#1}}
\newcommand{\bex}[1]{\begin{example}\label{#1}}
\newcommand{\be}[1]{\begin{equation}\label{#1}}
\newcommand{\ec}{\end{corollary}}
\newcommand{\efact}{\end{fact}}
\newcommand{\eex}{\end{example}}
\newcommand{\el}{\end{lemma}}
\newcommand{\er}{\end{remark}}
\newcommand{\ers}{\end{remarks}}
\newcommand{\et}{\end{theorem}}
\newcommand{\ed}{\end{definition}}
\newcommand{\ep}{\end{proposition}}
\newcommand{\ee}{\end{equation}}

\newcommand{\R}{\mathbb{R}}

\newcommand{\E}{\mathbb{E}}
\newcommand{\N}{\mathbb{N}}

\newcommand{\Z}{\mathbb{Z}}
\newcommand{\Q}{\mathbb{Q}}
\newcommand{\pp}{\mathbb{P}}
\newcommand{\kA}{\mathcal{A}}

\newcommand{\kC}{\mathcal{C}}
\newcommand{\kG}{\mathcal{G}}
\newcommand{\kR}{\mathcal{R}}
\newcommand{\kO}{\mathcal{O}}
\newcommand{\kP}{\mathcal{P}}

\newcommand{\kF}{\mathcal{F}}
\newcommand{\kE}{\mathcal{E}}
\newcommand{\kH}{\mathcal{H}}
\newcommand{\kI}{\mathcal{I}}
\newcommand{\kS}{\mathcal{S}}
\newcommand{\kN}{\mathcal{N}}
\newcommand{\kM}{\mathcal{M}}

\newcommand{\kK}{\mathcal{K}}

\newcommand{\lin}{\left[\kern-0.15em\left[}
\newcommand{\rin} {\right]\kern-0.15em\right]}
\newcommand{\ilin}{\left]\kern-0.15em\left]}
\newcommand{\irin} {\right[\kern-0.15em\right[}

\def\reff#1{(\ref{#1})}
\def \ind {\hbox{1\hskip -3pt I}}
\def\capa{\text{cap}}

\def\bs{\backslash}
\newcommand {\sous}[1] {\underline{#1}}

\begin{document}

\title{{\bf Moderate deviations for the range of a transient
random walk: path concentration}}

\author{
\normalsize{\textsc{Amine Asselah}\footnote{Aix-Marseille Universit\'e \& 
Universit\'e Paris-Est Cr\'eteil
E-mail: amine.asselah@u-pec.fr}\ \ \ \&\ \ \textsc{Bruno
Schapira}\footnote{Aix-Marseille Universit\'e, CNRS, Centrale Marseille, I2M, UMR 7373, 13453 Marseille, France. E-mail: bruno.schapira@univ-amu.fr}}}

\date{}
\maketitle

\begin{abstract}
We study downward deviations of the boundary of the range
of a transient walk on the Euclidean lattice. We describe
the optimal strategy adopted by the walk in order to
shrink the boundary of its range. The technics we develop
apply equally well to the range, and provide pathwise statements for
the {\it Swiss cheese} picture of Bolthausen, van den Berg and
den Hollander \cite{BBH}.
\newline
MSC 2010 \emph{subject classifications.} Primary 60F10; 60G50.
\end{abstract}

\section{Introduction}
In this paper we study downward deviations of the boundary of the range
of a simple random walk $(S_n,\ n\in \N)$ on $\Z^d$, with $d\ge 3$. 
The range at time $n$, denoted $\kR_n$, is the set of visited 
sites $\{S_0,\dots,S_n\}$, and
its boundary, denoted $\partial \kR_n$, is the set of sites of $\kR_n$
with at least one neighbor outside $\kR_n$. 
Our previous study \cite{AS} focused on the typical behavior
of the boundary of the range, whereas this work is devoted
to downward deviations and applications to a hydrophobic polymer model. 
The zest of the paper is about describing
the optimal strategy adopted in order to 
shrink the boundary of the range,
and our approach shed some light on the shape of the walk
realizing such a deviation.
In \cite{AS}, we emphasized the ways in which, for a transient
walk, the range and its boundary share a similar nature. 
Thus, even though the boundary of the range is our primary
interest, we mention at the outset that the technics we develop
apply equally well to the range. Since this last issue has been the
focus of many celebrated works, let us describe first the state
of the art there.

\paragraph{Deviations of the range.}
A pioneering large deviation study of Donsker and Varadhan
\cite{DV1} establishes asymptotics for downward deviations of the
volume of the Wiener sausage $t\mapsto W^a(t)$, 
that is the Lebesgue measure of an $a$-neighborhood of
the standard Brownian motion. The main result of \cite{DV1}
establishes,
in any dimension and for any $\beta>0$, the following asymptotics
\be{DV-1}
\lim_{t\to\infty}
t^{-\frac{d}{d+2}}\log\E[\exp\big(-\beta W^a(t)\big)]=
\frac{d+2}{2}\ \beta
(\frac{2\lambda_D}{d\beta})^{\frac{d}{d+2}},
\ee
where $\lambda_D$ is the first eigenvalue of the Laplacian
with Dirichlet condition on the boundary of a sphere
of volume one.
The asymptotics \reff{DV-1}, obtained in the random
walk setting in \cite{DV2}, correspond to downward deviation
of the volume of the range $\{|\kR_n|\le f(n)\}$ where
$|\kR_n|$ denotes the volume of $\kR_n$ and 
$f(n)$ is of order $n^{\frac{d}{d+2}}$. 
They suggest that during time $n$ a random walk is localized
in a ball of radius $(n/\beta)^{\frac{1}{d+2}}$ filled without holes. 
Bolthausen \cite{Bolt1} and Sznitman 
\cite{Sznitman1}, with different technics, extended the result of
\cite{DV1} to cover downward deviations corresponding to
$f(n)=n^{1-\delta}$ for any $\delta>0$.
A consequence of their analysis is that for $0<\gamma\le 2$
\be{Bolt-Sznit}
\lim_{t\to\infty}
t^{-\frac{d+\gamma-2}{d+\gamma}}
\log\E[\exp\big(-\beta\ t^{-\frac{2-\gamma}{d+\gamma}} W^a(t)
\big)]= -\frac{d+2}{2}\beta
(\frac{2\lambda_D }{d\beta})^{\frac{d}{d+2}}.
\ee
Then, three deep studies dealt with
the trajectory conditioned on realizing a large deviation
by Sznitman \cite{Sznitman2}, Bolthausen \cite{Bolt2} and
Povel \cite{Povel}. The case $\gamma=0$ in \reff{Bolt-Sznit}
was recognized as {\it critical} by
Bolthausen~\cite{Bolt1}, and indeed a different behavior
was later proved to hold \cite{BBH}.
The series of papers on downward deviations
culminated in a paper of Bolthausen, van den Berg and den Hollander
\cite{BBH} which covers the {\it critical regime}
$\{|\kR_n|-\E[|\kR_n|]\le \varepsilon n\}$. The latter contribution offers
a precise Large Deviation Principle, but no pathwise statement
characterizing the most likely scenario. The present paper is
a step towards filling this gap and providing answers
to their motto {\it How a Wiener sausage turns into
a Swiss cheese}? Let us quote their mathematical results.
In dimension $d\ge 3$, $\E[W^a(t)]$ grows linearly
and the limit of $\frac{1}{t}\E[W^a(t)]$ is denoted $\kappa_a$
(the Newtonian capacity of a ball of radius $a$).
It is proved in \cite{BBH} that for any $0<\varepsilon<1$ 
\be{BBH-1}
\lim_{t\to\infty}
\frac{1}{t^{\frac{d-2}{d}}}\log \pp\big[
W^a(t)-\E[W^a(t)]\le -\varepsilon \kappa_a t\big]=-I_a(\varepsilon),
\ee
where
\be{BBH-2}
I_a(\varepsilon)=
\frac{1}{2\kappa_a^{2/d}} \inf\{\|\nabla f\|_2:\ 
f\in H^1(\R^d), \|f\|_2=1,\ \int_{\R^d} 
\big(1-\exp(-f^2(x)))dx\ \le 1-\varepsilon\}.
\ee
A similar result for simple random walks is obtained in
Phetpradap's thesis \cite{Phet}: $\kappa_a$ becomes the
non-return probability say $\kappa_d$, and the factor $1/2\kappa_a^{2/d}$
in \reff{BBH-2} becomes $1/2d\kappa_d^{2/d}$.
When $d=3$ or $d=4$, the minimizers of \reff{BBH-2} are
strictly positive on $\R^d$, and decrease in the radial component.
This is interpreted as saying that Wiener sausage {\it "looks like
a Swiss cheese" with random holes whose sizes are of order 1 and
whose density varies on scale $t^{1/d}$}. On the other hand,
when $d\ge 5$, and when the parameter $\varepsilon$ in \reff{BBH-1}
is small, there is no minimizer for the variational problem
\reff{BBH-2}, suggesting that {\it the optimal strategy
is time-inhomogeneous}.

\paragraph{Boundary of the range.}
The boundary of the range, in spite of not receiving much attention,
enters naturally into the modelling of {\it hydrophobic polymers}.
Indeed, a {\it polymer} is a succession of monomers centered at the
positions of the walk (and thus covering $\kR_n$),
the complement of the range is occupied by the aqueous solvent,
and being {\it hydrophobic} means that
the monomers try to hide from it. A natural model is then 
the following polymer measure depending on two parameters: its length
$n$, and its inverse temperature $\beta$.
\[
d\widetilde \Q_n^\beta=\frac{1}{\widetilde Z_n(\beta)} 
\exp(-\beta |\partial\kR_n|) \, d\pp_n,
\]
where $\pp_n$ denotes the law of the simple random walk up to time $n$ and $\widetilde Z_n(\beta)$, the partition function, is a normalizing factor.
Biology suggests that as one tunes $\beta$, for a fixed polymer
length, a phase transition appears. The recent results of
Berestycki and Yadin \cite{BY} treat an asymptotic regime of length
going to infinity, and suggest that for any positive $\beta$ a long
enough polymer, that is under $\widetilde \Q_n^\beta$, is localized in a ball
of radius $\rho_n$ with $\rho_n^{d+1}$ of order $n$.
Thus, to capture the insight from Biology, we rather
scale $\beta$ with $n^{2/d}$, when $n$ is taken to infinity. 
We therefore consider
\be{def-Q}
d\Q_n^\beta=\frac{1}{Z_n(\beta)} \exp\Big(-
\frac{\beta}{n^{2/d}}\big(|\partial\kR_n|-\E[|\partial\kR_n|]\big)\Big)\, d\pp_n. 
\ee
The centering of $|\partial\kR_n|$ is a matter of taste, but the
scaling of $\beta$ by $n^{2/d}$ is crucial, and corresponds to
{\it a critical regime} for the boundary of the range reminiscent
of \reff{Bolt-Sznit} for $\gamma=0$. Indeed, understanding
the polymer measure is linked with analyzing the scenarii 
responsible for shrinking the boundary of the range on the scale
of its mean. However, before tackling deviations, let us
recall some typical behavior of the boundary of the range.
Okada \cite{Ok1} has proved a law of large numbers 
in dimension $d\ge 3$, and 
when dimension is two, he proved that  
\begin{equation*}
\frac{\pi^2}{2}\le \lim_{n\to\infty}
\frac{\E[|\partial\kR_n|]}{n/(\log n)^2} \le 2 \pi^2.
\end{equation*} 
Note that Benjamini, Kozma, Yadin and Yehudayoff \cite{BKYY}
in their study of the entropy of the range have obtained
the correct order of magnitude for $\E[|\partial\kR_n|]$ in $d=2$,
and have linked the entropy of the range to the size of its boundary.

In addition, a central limit theorem for the boundary
of the range was proved in \cite{AS} in dimension $d\ge 4$. 
When dimension is three, the variance is expected to grow
like $n\log n$, and only an upper bound of 
the right order is known \cite{AS}. 
We henceforth focus on the ways in which a random walk
reduces the boundary of its range. 

\paragraph{Capacity of the range.}
A key object used to probe the shape of the random walk is 
the {\it capacity of its range}. We first define it, and then state
our result.
For $\Lambda\subset \Z^d$, let $H_\Lambda^+$ be the time needed
by the walk to return to
$\Lambda$. The capacity of $\Lambda$, denoted $\capa(\Lambda)$, is
\be{def-capa}
\capa(\Lambda)=\sum_{x\in \Lambda} 
\pp_x\big[H_\Lambda^+=\infty\big].
\ee
Let us recall one of its basic property. There exists a positive constant $c_{\text{cap}}$, such that 
for all finite subset $\Lambda\subset \Z^d$ 
\be{cap.Lambda.general}
c_{\capa}\, |\Lambda |^{1-\frac 2d} \, \le\,  \capa(\Lambda)\, 
\le\,  |\Lambda|.
\ee
The upper bound follows by definition and the lower bound
is well known (see the proof of Proposition 2.5.1 in \cite{L}).
In a weak sense, the capacity of a set characterizes
the shape of a set: the closer it is to a ball, the smaller is 
its capacity. In view of \reff{cap.Lambda.general},
this is captured by the index $\kI_d$, defined  
for finite subsets $\Lambda$ in $\Z^d$, by 
\be{def-iso}
\kI_d(\Lambda):=\frac{\capa(\Lambda)}{|\Lambda|^{1-\frac 2d}}.
\ee
It is known that the index $\kI_d$ of a ball is bounded by some constant, independently of the radius of the ball (see \eqref{cap.boules} below).  
On the other hand, the capacity of the range $\kR_n$ has been studied in 
\cite{ABP, JO, L, RS}, and it is known that its mean is of order
$\sqrt n$ in dimension three, of order $n/\log n$ in
dimension four, and grows linearly in dimension five and larger.
Thus typically, for a transient walk,
$\kI_d(\kR_n)$ goes to infinity with $n$.

\paragraph{Results.}
We have found that the strategy for shrinking the (size of the)  
boundary of the range $|\partial \kR_n|$
below its mean by $\varepsilon n$ is different in dimension
three and in dimensions five and larger.
In dimension three (see \eqref{new-d3} below),
the walk has to spend a positive fraction of its time 
in a set with $\kI_d$-index of order 1 and volume of order $n/\varepsilon$, whereas 
in dimensions five and larger (see \eqref{new-d5} below), 
it has to spend a fraction of order $\varepsilon$ of its time in a set 
of $\kI_d$-index of order 1 and volume of order $n$.

In addition, 
we prove that spending a positive fraction (resp. a fraction $\varepsilon$) of the time in a ball of radius $(n/\varepsilon)^{1/d}$ (resp. $n^{1/d}$), leads to reducing the boundary of the range by a factor $\varepsilon$ in dimension three (resp. four and larger), and this gives us the lower bounds for the large deviations in \eqref{shrink-d3}, \eqref{shrink-d5} and \eqref{shrink-d4}, see also Section \ref{sec-low} for more details.

For a subset $\Lambda$
of $\Z^d$, we denote by $\ell_n(\Lambda)$ the
time spent by the walk inside $\Lambda$ up to time $n$.
For an integrable random variable $X$, we also denote by $\overline{X}$ 
the centered variable $X-\E[X]$. 

\bt{theo-shrink3}
Assume that $d=3$. There exist 
$\alpha\in (0,1)$ and $C>0$, such that for 
all $\varepsilon \in(0,\nu_3/2)$, 
\be{new-d3}
\lim_{n\to\infty} \pp\Big[
\exists \Lambda \subset \Z^3\, :\, \ell_n(\Lambda)\ge \alpha n,\,   
\kI_d(\Lambda)\le C, \, \frac{n}{C\varepsilon}\le
|\Lambda|\le C\frac n\varepsilon \ \Big|\   
\overline{|\partial\kR_n|}\le -\varepsilon n \Big]=1.
\ee
Moreover, there exist positive constants $\sous \kappa_3$ 
and $\bar\kappa_3$, 
such that for all $\varepsilon\in(0,\nu_3/2)$, and $n$ large enough
\be{shrink-d3}
\exp\big(-\bar \kappa_3\cdot \varepsilon^{\frac 23} n^{\frac 13}\big)\, 
\le\, 
\pp\Big[|\partial\kR_n|-\E[|\partial\kR_n|]\le -\varepsilon n\Big]\, \le\, 
\exp\big(-\sous\kappa_3\cdot \varepsilon^{\frac 23} n^{\frac13}\big).
\ee
\et

\bt{theo-shrink5}
Assume that $d\ge  5$. There exist 
$\alpha\in (0,1)$ and $C>0$, such that for all $\varepsilon\in (0,\nu_d/2)$,
\be{new-d5}
\lim_{n\to\infty} \pp\Big[
\exists \Lambda \subset \Z^d\, :\, \ell_n(\Lambda)\ge 
\alpha \varepsilon n,\, 
\kI_d(\Lambda)\le C,\, \frac{\varepsilon n}{C}\le |\Lambda|\le C\varepsilon n \ \Big|\   
\overline{|\partial\kR_n|}\le -\varepsilon n \Big]=1.
\ee
Moreover, there exist positive
constants $\sous \kappa_d$ and $\bar \kappa_d$, such that 
for all $\varepsilon \in (0,\nu_d/2)$, and $n$ large enough
\be{shrink-d5}
\exp\big(-\bar \kappa_d \cdot (\varepsilon n)^{1-\frac 2d}\big)\, \le\, 
\pp\Big[|\partial\kR_n|-\E[|\partial\kR_n|]\le -\varepsilon n\Big]\, \le\, 
\exp\big(-\sous \kappa_d\cdot (\varepsilon n)^{1-\frac 2d}\big).
\ee
\et

In dimension $4$, our result is slightly 
less precise: our upper bound for the large deviations 
is larger than the lower bound by a factor $|\log \varepsilon |^{1/2}$,  
and we only describe the localization phenomenon 
for $\varepsilon$ away from $0$ and $\nu_4$.
\bt{theo-shrink4}
Assume that $d=4$. 
For any $\varepsilon \in (0,\nu_4/2)$, 
there exist $\alpha=\alpha(\varepsilon) \in (0,1)$, 
and $C=C(\varepsilon)>0$, such that  
\be{new-d4}
\lim_{n\to\infty} \pp\Big[
\exists \Lambda \subset \Z^4\, :\, \ell_n(\Lambda) \ge \alpha n,  \, 
\kI_d(\Lambda)\le C,\, \frac nC\le |\Lambda|\le Cn \ \Big|\   
\overline{|\partial\kR_n|}\le -\varepsilon n \Big]=1.
\ee
Moreover, there exist 
positive constants $\sous \kappa_4$ and $\bar \kappa_4$, 
such that for all $\varepsilon \in (0,\nu_4/2)$, and $n$ large enough
\be{shrink-d4}
\exp\left(-\bar \kappa_4\cdot (\varepsilon n)^{\frac 12}\right) \, \le\,  
\pp\left[|\partial \kR_n|-\E[|\partial\kR_n|]\le -\varepsilon n\right]
\, \le\,  \exp\left(-\sous \kappa_4\cdot \frac{(\varepsilon n)^{1/2}}{|\log \varepsilon |^{1/2}} 
\right).
\ee
\et

\br{rem-Range}
Our theorems play with two parameters $n$ and $\varepsilon$,
and consider the regime where $n$ is large and $\varepsilon$ is small.
They suggest that conditioned on having a small
boundary of the range, in dimension three a (large) fraction of
the walk is localized in a ball of volume $n/\varepsilon$,
whereas in dimension five or larger a length $\varepsilon n$
of the trajectory is localized in a ball of volume 
of order $\varepsilon n$. 
This is the time-inhomogeneous nature of the trajectory that was mentioned earlier.  
\er

\br{rem-Proof}
We state the results for $|\partial \kR_n|$, 
and they stand as well for $|\kR_n|$. In this case,
\eqref{new-d3}, \eqref{new-d5}, and \eqref{new-d4} 
are new, and complement the results of 
van den Berg, Bolthausen, and den Hollander \cite{BBH}.
Our proof is based on a recursive slicing of a strand of
random walk, as we explain in Section~\ref{sec-prel}, whereas
the proof of \cite{BBH} consists of many steps:
(i) first a compactification obtained by
wrapping the trajectory on a Torus of side $An^{1/d}$,
with $A$ large, (ii) fixing the position of the walk at times
multiple of $n^{2/d}$, the so-called {\it skeleton}, and 
using concentration to showing that
the range and {\it averaged range conditioned on the skeleton} are close 
and finally (iii) representing the averaged range 
conditioned on the skeleton as a continuous functional of
the pair empirical measure, and invoking Donsker Varadhan
Large Deviation Theory. Note that step (ii) uses that
the range is a union of ranges over periods of length
$n^{2/d}$, and thus a Lipshitz function of the collection of
smaller ranges.
The boundary of the range is neither a union of boundaries, nor
a monotonous function of time, and we had to follow a rougher
but more robust method (see below for a sketch of our approach).
\er
\br{rem-CLT}
It is noted in \cite{BBH} that there is an anomaly in the scaling
of the rate function for $\varepsilon$ close to zero, which does
not connect with the central limit theorem. We quote from
\cite{BBH} ``{\it The anomaly for $d\ge 3$ 
is somewhat surprising. It suggests that the central
limit behavior is controlled by the local fluctuation [...],
while the moderate and large deviations are controlled
by the global fluctuations.}" This is exactly what our
decomposition shows: the central limit theorem obtained
in \cite{AS}, and the downward deviations obtained here correspond
to two distinct parts of the slicing down of a random trajectory:
the {\it self-similar independent parts} 
contribute to the fluctuations, whereas
{\it the mutual-intersection parts} contribute to downward deviations.
We will come back on this later.  \er

The capacity of the range plays a central role in our results.
Our main technical contribution, interesting on its own,
is the following estimate which generalizes an inequality 
(1.8) of \cite{AC}, as well as Proposition 1.5 of \cite{AA}.
We bound the probability of multiple visits to non-overlapping spheres,
all of the same radius, and our bound
involves the capacity of the union of the spheres. Let us first
introduce handy notation. For $x\in \Z^d$ and $r>0$, we denote
with $B(x,r)$ the Euclidean ball of radius $r$ and center $x$, and  
for any subset $\kC\subset \Z^d$, we let
\[
B(\kC,r):=\bigcup_{x\in \kC} B(x,r).
\]
Additionally, for $r>0$, $\kA(r)$ is the
collection of finite sets of centers defining non-overlapping spheres of radius $2r$
\be{def-kA}
\kA(r):=
\big\{ \kC\subset \Z^d:\ |\kC|<\infty\ \textrm{ and }\ \|x-y\|\ge 4r\textrm{ for all }x\neq y\in \kC\big\}.
\ee
\bp{prop.cap}
Assume that $d\ge 3$.
There exist positive constants $C$ and $\kappa$,  
such that for any $t>0$, $r\ge 1$, $\kC\in \kA(r)$, and $n\ge 1$, 
one has
\be{main-tech}
\pp\Big[\ell_n(B(x,r))\ge t \ \textrm{ for all }x\in \kC\Big]
\le  C\, (|\kC|\, n)^{|\kC|}  \, 
\exp\left(-\kappa \cdot \frac{t|\kC|}{|B(\kC,r)|^{\frac 2d} }\cdot
\kI_d(B(\kC,r))\right).
\ee
\ep
This result is useful
when the combinatorial term in \eqref{main-tech} is innocuous, and using \reff{cap.Lambda.general} we see that this holds 
when for a constant $\delta$ small
enough
\be{cond-tech}
|B(\kC,r)|^{2/d}\log n \, \le\,  \delta t.
\ee
Note also that on the right hand side of \reff{main-tech}, 
$t|\kC|$ is a lower bound for the total time spent in $B(\kC,r)$, under the event of the left hand side. 
If this total time is comparable with the volume of $B(\kC,r)$,
then the estimate discriminates sets of comparable volume
with different indices $\kI_d$. This is exactly the kind of situation we will encounter here.
To be more precise we will use \reff{main-tech} in cases where both $t|\kC|$ and $|B(\kC,r)|$ are of order $n$. 
In these cases, \eqref{cond-tech} holds when $r^d\ge C n^{2/d} \log n$, for some large enough constant $C$, and \eqref{main-tech} becomes 
 \be{ineq-cap1}
\pp\Big[\ell_n(B(x,r))\ge t\ \textrm{ for all }x\in \kC\Big]\ 
\le\ C \exp\big(- \kappa' \cdot \kI_d(B(\kC,r))\cdot n^{1-2/d}\big), 
\ee
with $\kappa'$ some other positive constant. Now in our results \eqref{new-d3}, \eqref{new-d5} and \eqref{new-d4}, the set 
$\Lambda$ will be a certain union of balls as here, and thus \eqref{ineq-cap1} explains why its $\kI_d$-index should be bounded.  

\paragraph{Sketch of our approach (in $d=3$).}
Our approach treats separately the upper and lower
bounds of the downward deviation in \reff{shrink-d3}. 
The upper bound relies on
the inequality \reff{lower-slicing} which roughly
bounds the centered boundary
of the range by self-similar terms, whose deviations are costly,
minus an increasing process $\xi_n(T)$ (defined in \reff{def-xi})
measuring the {\it folding of the trajectory} over small strands of
length $T$. Since this term is
the key term driving the moderate deviations, 
let us mention that it can be written 
in terms of Green's function up to time $T$, denoted $G_T$ and defined
later (and the $\sim$ means that the {\it correct} definition of $\xi_n(T)$
requires more notation),
\[
\xi_n(T) \sim \frac 1T \sum_{k=1}^n \sum_{z\in \kR_k}  G_T(z-S_k).
\]
What we show first is that a decrease in $\partial \kR_n$ most
likely translates into an increase of $\xi_n(T)$, for an appropriate
scale $T$ . Then, we show that the event that 
$\xi_n(T)$ large means that {\it many}
balls are visited {\it often} by the walk: see \reff{no-random} 
and Lemma~\ref{lem-key} for a precise statement.
The probability of this latter event is controlled by the capacity
of the collection of balls as in Proposition~\ref{prop.cap}.
We believe that a similar structure shows up in the
capacity of the range in dimension five or more, which has been the focus 
of recent studies (see \cite{ABP} and references therein).
The lower bound has a different flavor, and relies on a
covering result saying that when the walk is localized in a ball of volume
$n/\varepsilon$, it likely visits an $\varepsilon$-fraction
of all fixed {\it large} subsets of the ball 
(see Proposition \ref{prop-cover}).

\paragraph{Application to the polymer measure.}
Recall that 
$$Z_n(\beta) = \E \left[ \exp\Big(-
\frac{\beta}{n^{2/d}}\big(|\partial\kR_n|-\E[|\partial\kR_n|]\big)\Big)\right].$$
Note first that by Jensen's inequality, one has $Z_n(\beta)\ge 1$, for all $\beta\ge 0$. 
Thus one can define for $\beta\ge 0$, 
$$F^+(\beta) = \limsup_{n\to \infty}\,  \frac{1}{n^{1-\frac 2d} }\, \log Z_n(\beta),$$
and 
$$
F^-(\beta) = \liminf_{n\to \infty}\,  \frac{1}{n^{1-\frac 2d} }\, 
\log Z_n(\beta).
$$
Denote by $\kA(r,v)$ (with $v>0$) 
the subset of $\kA(r)$, whose elements $\kC$ satisfy 
$|B(\kC,r)|\le v$.

\bt{theo-poly}
Assume that $d\ge 3$. The following statements hold. 
\begin{enumerate}
\item The functions $F^+$ and $F^-$ are non decreasing in $\beta$. As a consequence there exist $0\le \beta_d^+\le \beta_d^-\le +\infty$, such that 
$$F^\pm(\beta)=0\quad \text{for }\beta<\beta_d^\pm \quad \text{and} \quad F^\pm(\beta)>0 \quad \text{for }\beta>\beta_d^\pm.$$ 
\item In fact $\beta_d^+$ and $\beta_d^-$ are positive and finite.
\end{enumerate}
Let $(r_n,\ n\in \N)$ be any fixed sequence of reals 
satisfying $n^{2/d}(\log n)^2\le r_n^d\le \frac{n}{\log n}$, 
for all $n\ge 1$. Then,
\begin{enumerate}
\item[3.] For any $\beta <\beta_d^+$, 
$\alpha \in (0,1)$ and $A>0$,
$$
\lim_{n\to\infty}
\Q_n^\beta \left[\exists \kC\in \kA(r_n,An) \, :\, \ell_n(B(x,r_n))\ge \alpha \frac{n}{|\kC|}\ \textrm{for all }x\in \kC \right] =0.
$$
Moreover, for any $\varepsilon \in (0,\nu_d)$, 
$$
\lim_{n\to\infty}
\Q_n^\beta \left[|\partial \kR_n|-\E[|\partial \kR_n|]
\le -\varepsilon n \right] = 0. 
$$
\item[4.] For any $\beta>\beta_d^-$, there exists $\varepsilon(\beta) \in (0,\nu_d)$, such that 
$$
\lim_{n\to\infty}
\Q_n^\beta \left[|\partial \kR_n|-\E[|\partial \kR_n|]
\le -\varepsilon(\beta) n \right] = 1. 
$$
Moreover, when $d=3$, there exist $\alpha\in (0,1)$ and $A>0$, such that
\be{poly-3}
\lim_{n\to\infty}
\Q_n^\beta \left[ \exists\kC\in \kA(r_n,\frac{An}{\varepsilon(\beta)})\, :\,  
\kI_d(B(\kC,r_n))\le A \text{ and }
\ell_n(B(x,r))\ge \alpha \frac{n}{|\kC|} \ \textrm{for all }x\in \kC  \right] = 1.
\ee
When $d\ge 4$, there exist $\alpha\in (0,1)$ and $A>0$, such that 
\be{poly-4}
\lim_{n\to\infty}
\Q_n^\beta \left[ \exists\kC\in \kA(r_n,An),\ \kI_d(B(\kC,r_n))\le A \text{ and }\ell_n(B(x,r_n))\ge \alpha\varepsilon(\beta)
\frac{n}{|\kC|} \ \textrm{for all }x\in \kC  \right] = 1.
\ee
\item[5.] In the previous part, one can choose $\varepsilon(\beta)$
such that $\varepsilon(\beta)\to \nu_d$, as $\beta \to +\infty$.  
\end{enumerate}
\et

The rest of the paper is organized as follows. Section~\ref{sec-prel}
sets the preliminaries, where Lemma~\ref{lem-key} is our key technical
tool, interesting on its own.
Section~\ref{sec-up} presents the upper bounds for our 
large deviation estimates, and Section~\ref{sec-low} presents
the corresponding lower bounds. A large part of Section~\ref{sec-low}
deals with the three dimensional case and proves 
Proposition~\ref{prop-cover} which establishes lower bound on
the probability of covering fraction of domains.
Section~\ref{sec-cap} introduces
capacities, and contains the proof of Proposition~\ref{prop.cap}, 
as well as establishing statements
\eqref{new-d3}, \eqref{new-d5} and \eqref{new-d4}.
Finally, our application to hydrophobic polymer, that is 
Theorem~\ref{theo-poly}, is explained in Section~\ref{sec-poly}.

\section{Preliminaries}\label{sec-prel}
\subsection{Notation}
For any integers $n\le  m$ we write 
$$\kR(n,m)=\left\{S_n,\dots,S_m\right\}.$$
We denote by $\|\cdot\|$ the Euclidean norm on $\R^d$, and by 
$$B(z,r)=\left\{y\in \Z^d\ :\ \|z-y\|\le r\right\},$$
the (discrete) ball of radius $r$ centered at $z$. 
For $z\in \Z^d$ and $\Lambda$ a subset of $\Z^d$, we write   
$$d(z,\Lambda)=\inf\{\|z-y\|\ :\ y\in \Lambda\},$$
for the distance between $z$ and $\Lambda$, and let $|\Lambda|$ be the size (which we also call the volume) of $\Lambda$. 
We also define  
$$\partial \Lambda = \left\{z\in \Lambda\ :\  d(z,\Lambda^c)=1\right\} \quad \text{and}\quad 
\Lambda^+= \Lambda\cup \left\{z\in \Z^d\ : \ d(z,\Lambda)=1\right\}.$$ 
Furthermore,  
$$H_\Lambda=\inf\{k\ge 0 \ :\ S_k\in \Lambda\},$$
is the hitting time of $\Lambda$, that we abbreviate in $H_z$ when $\Lambda=\{z\}$. and we recall that for $n\ge 0$, 
$$\ell_n(\Lambda)=\,  |\{k\le n\ :\ S_k\in \Lambda\}|,$$
is the number of steps spent in $\Lambda$ before time $n$. 
We next consider Green's function, 
$$G(z)=\sum_{k\ge 0} \pp[S_k=z]=\pp[H_z<\infty] \times G(0)\qquad \text{for }z\in \Z^d,$$
and its restricted version for any positive integer $T$
$$G_T(z)=\sum_{k=0}^T   \pp[S_k=z].$$

\subsection{Range and its boundary.}
\label{subsec.RBR}
A powerful idea, going back at least to Le Gall \cite{LG},
is to cut a trajectory into two pieces, and
read the total range as a union of two ranges minus the mutual 
intersection of two {\it independent} strands. From a set-theoretical
point of view, the volumes satisfy an exclusion-inclusion formula:
\be{heur-1}
|\Lambda_1\cup \Lambda_2|=|\Lambda_1|+|\Lambda_2|-|\Lambda_1\cap\Lambda_2| \qquad \text{for all }\Lambda_1,\, \Lambda_2\subset \Z^d. 
\ee
The key {\it probabilistic} point is that if $\Lambda_1=\{S_0,\dots,
S_n\}$ and $\Lambda_2=\{S_n,\dots,S_{n+m}\}$, then by translating
sets by $S_n$, we have
\be{heur-2}
|\Lambda_1\cap\Lambda_2|=\big|\{S_0-S_n,S_1-S_n,\dots,0\}\cap
\{0,\dots,S_{n+m}-S_n\}\big|,
\ee
and the two sets on the right hand side are independent.
It is then possible, using the symmetry of the walk, to compute the expectation:
\be{heur-exp}
\E[|\Lambda_1\cap\Lambda_2|]=\sum_{z\in \Z^d} \pp[H_z\le n]\, \pp[H_z\le m].
\ee
In view of \reff{heur-1} and \reff{heur-2}, it is no surprise that
Bolthausen, van den Berg and den Hollander 
studied at the same time downward deviations of the
Wiener sausage \cite{BBH} and upward deviation of the intersection
of two Wiener sausages in \cite{BBH2}. The second paper is
independently motivated by an older paper of Khanin, Mazel, Shlosman,
and Sinai \cite{KMSS} studying bounds for the intersection of
two independent ranges in an infinite time-horizon, a problem which is 
still open.

Now, the boundary of the range does not quite satisfy the
exclusion-inclusion equality \reff{heur-1}. However, 
if $z\in \partial \Lambda_1\bs \Lambda^+_2$, then the neighbors
of $z$ are not in $\Lambda_2$, and at least one of them is not in 
$\Lambda_1$. This means that $z\in \partial(\Lambda_1\cup\Lambda_2)$,
and therefore 
\be{heur-3}
\big( \partial \Lambda_1\bs \Lambda^+_2\big)\cup
\big( \partial \Lambda_2\bs \Lambda^+_1\big)\subset 
\partial\big( \Lambda_1\cup\Lambda_2\big).
\ee
When taking the volume, this reads
\be{heur-4}
\begin{split}
|\partial\big( \Lambda_1\cup\Lambda_2\big)|\ge &
|\partial\Lambda_1|+|\partial\Lambda_2|-
\Big(|\partial\Lambda_1\cap \Lambda^+_2|+
|\Lambda^+_1\cap\partial\Lambda_2|\Big)\\
\ge & |\partial\Lambda_1|+|\partial\Lambda_2|-
2|\Lambda^+_1\cap\Lambda^+_2|.
\end{split}
\ee
On the other hand, if $z\not\in \partial\Lambda_1$,
but $z\in \Lambda_1$, then $\{z\}^+\subset \Lambda_1$, and
$z$ is not in the boundary of $\Lambda_1\cup\Lambda_2$. In
other words,
\be{heur-5}
|\partial\big( \Lambda_1\cup\Lambda_2\big)|\le
|\partial\Lambda_1|+|\partial\Lambda_2|-|\Lambda_1\cap
\partial \Lambda_2|.
\ee
Thus, \reff{heur-4} and \reff{heur-5} make up for the equality
\reff{heur-1}.

By exploiting \eqref{heur-4}, it has been shown in \cite{AS} that there exist positive 
dimension-dependent constants $\nu_d$ and $C_d$, such that  
\be{mean-transient}
\big| \E[|\partial \kR_n|]-\nu_d n\big|\le C_d\, \psi_d(n),
\ee
for all $n\ge 1$, where 
\be{psid}
\psi_3(n)=\sqrt n,\qquad \psi_4(n)= \log n,\qquad \psi_d(n)= 1 \ \textrm{ for $d\ge 5$}.
\ee

\subsection{On Green's function}
First, we recall the asymptotics of
Green's function (see \cite[Theorem 4.3.1]{LL}) 
\be{Green}
G(z)=\kO\big(\frac{1}{1+\|z\|^{d-2}}\big). 
\ee
For Green's function restricted to the first $T$ steps, 
the following holds.
\bl{Green.restricted}
Assume that $d\ge 3$. There exist positive constants $c$ and $C$, such that for any $T>0$ and $z\in \Z^d$, 
$$G_T(z)\ \le\  C\,  \frac{T}{1+\|z\|^d}\, \exp\left(-c\, \frac{\|z\|^2}{T}\right).$$
\el
\begin{proof} One can assume that $\|z\|\ge \sqrt T$, 
as otherwise the result follows from \eqref{Green}.
One result of \cite{HSC} ensures that there exist 
constants $c$ and $C$, such that
\begin{eqnarray}
\label{kernel}
\pp(S_n=z) \le C\, \frac 1{n^{d/2}}\, \exp(-c\|z\|^2/n)\quad 
\textrm{for all $z$ and $n\ge 1$}.
\end{eqnarray}
This implies the following bound.
\begin{equation*}
\begin{split}
G_T(z) &=  \sum_{k=1}^T \pp[S_k=z]
\le\ C\sum_{k=1}^T k^{-d/2} \exp(-c \|z\|^2/k) \\ 
&\le\ C\, \|z\|^{-d} \sum_{k=1}^T \exp(-\frac c2 \|z\|^2/k)
\le\ C\, \frac{T}{ \|z\|^d}\, \exp(-\frac c2 \|z\|^2/T).
\end{split}
\end{equation*}
\end{proof}

For $x\in \Z^d$ and $r\in \N$, we consider the discrete cube
centered on $x$ and of side $r$:
\[
Q(x,r)=\big(x+]-r,r]^d\big)\cap \Z^d.
\]
The elements of a partition of $\Z^d$ obtained from translates
of $Q(0,r)$ is denoted $\kP_r$ with
\be{def-A}
\kP_r=\{Q(x,r):\ x\in 2r\Z^d\}.
\ee
For $T>0$, we denote by $\kP_r(T)$ the
elements of $\kP_r$ whose intersection with $B(0,T)$ is not empty. 
We will need the following covering result.
\bl{lem-green}
Assume that $d\ge 3$. 
There exists a constant $C>0$, such that for any $r\ge 1$, $\gamma \in (0,1)$, and any collection of subsets $\{\Lambda_Q,\, Q\in \kP_r\}$ satisfying 
\be{hyp-lem.green}
\Lambda_Q\subset Q\quad
\text{and}\quad |\Lambda_Q|\le \gamma \, r^d\qquad 
\text{for all }Q\in \kP_r,
\ee
one has for all $T\ge r$, 
\be{green-main}
\sum_{Q\in \kP_r} \sum_{z\in \Lambda_Q} G_T(z) \ 
\le\  C\,  \big( \gamma^{2/d} r^2+\gamma T\big).
\ee
\el
Note that since the elements of $\kP_r$ form a partition of the space (in particular they are disjoints), 
the left hand side of \eqref{green-main} is also equal to the expectation of the time spent before $T$ in the union of all the $\Lambda_Q$, with  $Q\in \kP_r$. In particular, when $\Lambda_Q=\Lambda\cap Q$ for some set $\Lambda$, then the left hand side of \eqref{green-main} is simply the expected value of the number of steps spent in $\Lambda$ before time $T$.  

\begin{proof}
We denote by $C$ a constant whose value might change from line to line.  
We decompose the sum in the left hand side of \reff{green-main} into
$\Sigma_I+\Sigma_{I\!I}$, with 
\begin{equation*}
\Sigma_I:=\sum_{Q\in \kP_r(\sqrt T)} \sum_{z\in \Lambda_Q}G_T(z)\quad\text{and}\quad
\Sigma_{I\!I}:=\sum_{Q\in \kP_r(T)\setminus \kP_r(\sqrt T)}\sum_{z\in \Lambda_Q}G_T(z).
\end{equation*}
For $\Sigma_I$ we use \eqref{Green} and \eqref{hyp-lem.green}. This gives 
\begin{equation*}
\begin{split}
\Sigma_I\le & \sum_{z\in \Lambda_{Q(0,r)}}\! G_T(z)+
\sum_{Q\in \kP_r(\sqrt T)\setminus \{Q(0,r)\}}\sum_{z\in \Lambda_Q} G_T(z)\\
\le & \ C\left( |\Lambda_{Q(0,r)}|^{2/d}+ \sum_{Q\in \kP_r(\sqrt T)\setminus\{Q(0,r)\}}
\frac{\gamma r^d}{d(0,Q)^{d-2}}\right)\\
\le & \ C \left(\gamma^{2/d}r^2+ \gamma r^d
\sum_{k=1}^{[\sqrt T/r]+1}\frac{k^{d-1}}{(kr)^{d-2}}\right)\\
\le &\  C\big(\gamma^{2/d}r^2+\gamma T\big),
\end{split}
\end{equation*}
where at the second line we used a well known bound for the first sum (see for instance the proof of Proposition 2.5.1 in \cite{L}). 
We deal now with $\Sigma_{I\!I}$ and use this time Lemma \ref{Green.restricted} instead of \eqref{Green}. This gives
\begin{equation*}
\begin{split}
\Sigma_{I\!I}\le &\  C\,   \gamma r^d T \sum_{Q\in \kP_r(T)\setminus \kP_r(\sqrt T)} 
\frac{\exp\big(-cd(0,Q)^2/T\big)}
{d(0,Q)^d}\\
\le & \ C\,  \gamma r^d T \sum_{k=[\sqrt T/r]}^{[T/r]+1} k^{d-1}
\frac{\exp(-ck^2r^2/T) }{(kr)^d}\\
\le & \ C\,  \gamma  T  \sum_{k\ge [\sqrt T/r]} \frac 1k \, \exp\left(-c\frac{k^2r^2}{T}\right) \\
 \le &\   C\,   \gamma T.
\end{split}
\end{equation*}
\end{proof}

\subsection{Rolling a Ball}
The main result of this section requires further notation. 
For $V$ a subset of $\Z^d$, define 
\be{def.K}
\kK_n(V,t):=\{k\in \{1,\dots,n\}:\ \ell_k(S_k+V)>t\}.
\ee 
Thus, as the set $V$ rolls along the random walk trajectory,
$\kK_n(V,t)$ records the times before $n$ 
when the number of visits in the  
moving window $V$ exceeds $t$. We use here only
the case where $V$ is a ball, and our first concern
is to express $\{|\kK_n(V,t)|>L\}$ in terms of occupation
of non-overlapping balls.
\bl{lem-prep}
For any $v\in \Z^d$, $r\ge 1$, $t>0$, $L>0$ and $n\ge 1$, 
there is a random subset
$\kC\in \kA(r)$,
which is measurable with respect to $\sigma(S_1,\dots,S_n)$, such that
\be{inclusion-prep}
\{|\kK_n(B(v,r),t)|>L\}\ \subset\ 
\{\ell_n(B(x+v,r)\ge t\ \textrm{ for all }  x\in \kC\} \, \cap\, 
\{\ell_n\big( B(\kC,4r)\big)\ge L\}.
\ee
\el
\begin{proof}
Assume $\{|\kK_n(B(v,r),t)|>L\}$ and define 
\[
k_1:=\inf\{k\ge 0\ :\ \ell_k(S_k+B(v,r))>t\}.
\]
Since $L>0$, we have $k_1\le n$, and we set $\kC_1=\{S_{k_1}\}$.
If all the elements $S_k$ with $k\in \kK_n(B(v,r),t)$ are in $B(\kC_1,4r)$, 
set $\kC=\kC_1$. Otherwise, proceed by induction and assume that for $i>1$,
$\kC_{i-1}\in \kA(r)$ and $k_{i-1}$ satisfy
\be{hypn.t}
\ell_{k_{i-1}}\Big(B(x+v,r)\Big)
>t \qquad \text{for all }x\in \kC_{i-1}.
\ee
If $\ell_n(B(\kC_{i-1},4r))\ge L$, then $\kC=\kC_{i-1}$. 
Otherwise, define 
$$
k_i:=\inf\{k>k_{i-1}:\ S_k\not\in B(\kC_{i-1},4r) 
\text{ and } \ell_k(S_k+B(v,r))>t\}.
$$
Note that $k_i\le n$, and define $\kC_i=\kC_{i-1}\cup \{S_{k_i}\}$. Since
$k_i\ge k_{i-1}+1$, the construction stops in a finite number of steps 
and yields a finite set $\kC\in \kA(r)$.
\end{proof}

We introduce now some new notation. For any positive integers $n$ and $m$, and any positive 
$r$, $t$ and $L$,  define 
\be{def-Gn}
\kG_n(r,t,m):= \bigcup_{\kC\in \kA(r)\, :\, |\kC|=m} \  
\left\{ \ell_n\Big(B(x,r)\Big) > t\ \textrm{ for all }x\in \kC\right\},
\ee
and  
\be{def-Hn}
\kH_n(r,L,m):=\bigcup_{\kC\in \kA(r)\, :\,  |\kC|=m}
\left\{\ell_n\Big(B(\kC,4r)\Big)\ge L\right\}.
\ee
Note that $\kG_n(r,t,m)\subset \kG_n(r,t,m-1)$ whereas
$\kH_n(r,L,m-1) \subset \kH_n(r,L,m)$.
Therefore \reff{inclusion-prep} implies that for {\it any} positive integer $m$
\be{no-random}
\{|\kK_n(B(v,r),t)|>L\}\ \subset\ \kG_n(r,t,m)\cup \kH_n(r,L,m). 
\ee
We recall next a useful inequality proved in \cite[Lemma 1.2]{AC}.
There exist positive constants $\kappa_0$ and $C$, such that  
for all $t>0$ and all (nonempty) subsets $\Lambda\subset \Z^d$,  
\be{Aminebound}
\pp\left[\ell_\infty(\Lambda) \ge t \right]\, 
\le\,  C\, \exp\left(-\kappa_0\cdot  \frac{t}{|\Lambda|^{2/d}}\right).
\ee
We are now ready to estimate the probability of the event $\{|\kK_n(B(v,r),t)|>L\}$.

\bl{lem-key}
Assume that $d\ge 3$. 
There exist positive constants $\kappa$ and $C$, such that for any 
$r\ge 1$, $L\ge 1$, $t\ge 1$, and $n\ge 2$, satisfying
\be{cond-key}
\left(1+\frac Lt\right)^{2/d} \log n\, \le\,  \kappa\,  \frac{t}{r^2}, 
\ee
we have for any $v\in \Z^d$, 
\be{ineq-key}
\pp\big[|\kK_n(B(v,r),t)|\ge L\big] \le \, C\, 
\exp\left\{-\kappa\,  \frac{t}{r^2}
\left(1+\frac{L}{t}\right)^{1-\frac{2}{d}}\right\}.
\ee
Moreover, for all $K>0$, there exists $\delta\in (0,1)$, such that
\be{KnGn}
\pp\Big[ |\kK_n(B(v,r),t)|\ge L, \, 
\kG^c_n(r,t,[\delta L/t])\Big] \, 
\le\,  C\, \exp\left\{-K\, \frac{t}{r^2}
\left(1+\frac{L}{t}\right)^{1-\frac{2}{d}}\right\},
\ee
with the convention that $\kG_n^c(r,t,0)$ is the empty set. 
\el
\begin{proof}
We start with the proof of \eqref{ineq-key}. 
First note that if $0<L<t$, then on the event $\{|\kK_n(B(v,r),t)|\ge L\}$, 
there is $k\le n$ such that $ \ell_n(B(S_k+v,r))>t$,
and then \reff{Aminebound} gives the result. Thus we can assume that $L\ge t$, and in fact also that $L\le n$ (as otherwise there is nothing to prove). Similarly by taking $\kappa$ small enough, one can assume that $r<n$.  
Choose now $m^*=1+[L/t]$ and note that if $\kappa$ is small enough, then \eqref{cond-key} implies that 
\[
(m^*)^{2/d}\log(4n)\, \le\,  \frac{\kappa_0}{18d} \, \frac{t}{r^2},
\]
with $\kappa_0$ as in \eqref{Aminebound}.
Now, inclusion \reff{no-random} gives 
\be{add-step}
\pp\big[|\kK_n(B(v,r),t)|>L\big] \le
\pp\big[\kH_n(r,L,m^*)\big]+ \pp\big[\kG_n(r,t,m^*)\big].
\ee
By \eqref{Aminebound},
(noting that $B(x,r)$ can be visited 
before time $n$ only if $\|x\| < 2n$, and that its cardinality is smaller than $(3r)^d$), 
\be{borneGn}
\begin{split}
\pp[\kG_n(r,t,m^*)] \le\ & \sum_{\kC\in \kA(r)\, :\,  |\kC|=m^*}
\pp\left[\ell_n\Big(B(\kC,r)\Big)\ge m^*t \right] \\
 \le \ &  C(4n)^{dm^*} \exp\left(-\frac{\kappa_0}{9}\, 
(m^*)^{1-\frac 2d}\, \frac{t}{r^2} \right) \\
 \le \ &  C\exp\left(-\frac{\kappa_0}{18}\, (m^*)^{1-\frac 2d}\, \frac{t}{r^2} \right).
\end{split}
\ee 
Likewise, if \eqref{cond-key} holds with $\kappa$ small enough, one has 
$$(m^*)^{1+2/d} \log (10n)\, \le\,  \frac{\kappa_0}{18\cdot 16d} \, \frac{L}{r^2},$$
and as in \eqref{borneGn} this leads, for any $m\le m^*$, to
\be{borneHn}
\begin{split}
\pp[\kH_n(r,L,m)] \ & \le\  \sum_{\kC\in \kA(r)\, :\,  |\kC|=m}
\pp\left[\ell_n\Big(B(\kC,4r)\Big)\ge L\right] \\
 &\le \   C\, (10n)^{dm}\,  \exp\left(-\frac{\kappa_0}{9\cdot 16}\, 
\frac{L}{m^{2/d}r^2} \right) \\
 &\le \   C\, \exp\left(-\frac{\kappa_0}{18\cdot 16}\, 
\frac{L}{m^{2/d}r^2} \right)\\
 &\le \  C\,  \exp\left(-\frac{\kappa_0}{36\cdot 16}\cdot 
\frac{1+L/t}{m^{2/d}}\cdot \frac{t}{r^2} \right),
\end{split}
\ee
where at the last line we used that $L/t\ge (1+L/t)/2$ (which holds since $L/t\ge 1$). 
Then \eqref{ineq-key} follows from \eqref{add-step}, 
\eqref{borneGn} and \eqref{borneHn}.

To prove \eqref{KnGn}, note that \eqref{no-random} implies that
\[
\{ |\kK_n(B(v,r),t)|\ge L\}\cap
\kG^c_n(r,t,[\delta L/t])\ \subset \ \kH_n(r,L,[\delta L/t]).
\]
Then take $\delta$ small enough,
so that $\kappa_0/(36\cdot 16\delta^{2/d})>K$, and use \eqref{borneHn} with $m=[\delta L/t]$. 
\end{proof}

\section{Upper Bounds}\label{sec-up}
\subsection{Slicing of a trajectory and first estimates}
The main idea is that by slicing a trajectory into
small pieces, the boundary of the range divides into a sum
of boundaries of these pieces minus mutual intersections.
The main result here is the following lower bound on $|\partial \kR_n|$.
\bp{prop-slicing}
For any positive integers $T$ and $n$ with $T\le n$, we have
\be{lower-slicing}
\overline{|\partial \kR_n|}\ \ge\ \kS_n(T)+\kM_n(T) -\frac{4d}{T}\sum_{k=1}^n
\sum_{z\in \kR_k^{++}} G_T(z-S_k)+\kE_n(T),
\ee
where $\kS_n(T)$ behaves like a sum of $n/T$ independent centered terms
bounded by $T$, and $\kM_n(T)$ behaves like a martingale with increments bounded by $T$. In particular   
there exists a positive constant $c$, such that for any positive $\varepsilon$, any $T=T(n)$ going to infinity with $n$, and $n$ large enough, we have
\be{prop-S}
\pp\big[\kS_n(T)+\kM_n(T)\le -\varepsilon n\big]\le
T\exp\big(- c\frac{\varepsilon^2 n}{T}\big).
\ee
Moreover, $\kE_n(T)=\kO(T+n\psi_d(T)/T)$, with $\psi_d(\cdot)$ as in \eqref{psid}.
\ep

\begin{proof}
Recall that for any subsets $\Lambda_1$ and $\Lambda_2$ in $\Z^d$, one has 
$$|\partial (\Lambda_1\cup \Lambda_2)| \ge |\partial \Lambda_1|+ |\partial \Lambda_2| - 2|\Lambda_1^+\cap \Lambda_2^+|.$$
Now, if we consider $N$ subsets $\{\Lambda_1,\dots,\Lambda_N\}$,
by induction one obtains 
\be{slice-3}
|\partial\big(\Lambda_1\cup\dots\cup\Lambda_N\big)|\ge
\sum_{i=1}^N |\partial \Lambda_i|-2\sum_{i=2}^N|
\Lambda_i^+\cap\big(\cup_{j<i} \Lambda_j^+\big)|.
\ee
Applying this to the boundary of the range, this gives 
for any positive integers $T \le n$, 
any $i\in \{-1,\dots,T-2\}$, and with the notation 
$K_n(T)=\lfloor n/T\rfloor -2$, 
\be{key.min}
|\partial \kR_n| \ge \left(\sum_{j=0}^{K_n(T)} 
|\partial \kR(i+jT+1,i+(j+1)T)|\right)  - X_{K_n(T)}(i,T) -\kO(T),
\ee
with for $k\le {K_n(T)}$
\be{def.Xn}
X_k(i,T)= 2\sum_{j=1}^k
|\kR^+_{i+jT}\cap \kR^+(i+jT+1,i+(j+1)T)|.
\ee
Moreover, the elements of the sum in \eqref{key.min} are i.i.d. and distributed like $|\partial \kR_{T-1}|$. 
For simplicity, denote the $j$-th term 
of this sum by $U_{i,j}(T)$. 
Now, using \eqref{mean-transient} leads to 
\be{key.min2}
\overline{|\partial \kR_n|} \ge \left(\sum_{j=0}^{K_n(T)}
\overline U_{i,j}(T) \right)  - X_{K_n(T)}(i,T) -
\kO\left(T+ \frac{n}{T}\psi_d(T)\right). 
\ee
Denote by $\kF_k(i,T)$ the $\sigma$-field generated by
$\{S_0,\dots,S_{i+kT}\}$ for $1\le k\le K_n(T)$.
Then, for $k\le K_n(T)$, define $M_k(i,T)$ by 
\be{MniT}
X_k(i,T)= M_k(i,T) +2\sum_{j=1}^{k}
\E\left[ |\kR^+_{i+jT}\cap \kR^+(i+jT+1,i+(j+1)T)| \, 
\Big|\, \kF_j(i,T)\right].
\ee
Note that for each $(i,T)$, $(M_k(i,T),\ k\le K_n(T))$ is an
$(\kF_k(i,T), \ k\le K_n(T))$-martingale. 
On the other hand, using that 
$$|\Lambda_1^+\cap \Lambda_2^+| \le 2d\, |\Lambda_1^{++}\cap \Lambda_2| \qquad \text{for all nonempty }\Lambda_1, \Lambda_2\subset \Z^d,$$ 
with $\Lambda_1^{++}=(\Lambda_1^+)^+$, we get 
\be{RniT}
\begin{split}
 \E\left[ |\kR^+_{i+jT}\cap \kR^+(i+jT+1,i+(j+1)T)| \, 
\Big|\, \kF_{i+jT}\right] \le  & \  2d  
\sum_{z\in \kR_{i+jT}^{++}} \pp_{S_{i+jT}}
\big[1\le H_z\le T\big]\\
\le & \ 2d \sum_{z\in \kR^{++}_{i+jT}} G_T(z-S_{i+jT}).
\end{split}
\ee
Now, define 
$$
\kS_n(T) = \frac{1}{T}\sum_{i=-1}^{T-2}
\sum_{j=0}^{K_n(T)} \overline U_{i,j}(T)
\quad \text{and}\quad
\kM_n(T) = - \frac 1T \sum_{i=-1}^{T-2} M_n(i,T),$$
so that using \eqref{key.min2}, \eqref{MniT}, and \eqref{RniT}, 
we establish \reff{lower-slicing}.
Now assume that $T=T(n)$ satisfies $T(n)\to\infty$, as $n\to \infty$. 
Then, using that $|U_{i,j}(T)|\le T$, and Azuma's inequality, we see that for any $i$, $T$, 
and $n$ large enough
\be{Uij}
\pp\left[\sum_{j=0}^{K_n(T)} \overline U_{i,j}(T)\le -
\frac{\varepsilon n}{2}\right] \le \exp\left(-c 
\frac{\varepsilon^2 n}{T}\right),
\ee
for some constant $c>0$. 
Therefore a union bound gives 
\be{Uij2}
\pp\left[\frac 1T \sum_{i=-1}^{T-2}\sum_{j=0}^{K_n(T)} 
\overline U_{i,j}(T)\le 
-\frac{\varepsilon n}{2}\right] \le T\exp\left(- c\frac{\varepsilon^2 n}{T}\right).
\ee
Likewise, one obtains 
\be{MnT}
\pp\left[M_n(T)\ge \frac{\varepsilon n}{2}\right]
\le T\exp\left(-c\, \frac{\varepsilon^2 n}{T}\right).
\ee
Inequality \reff{prop-S} follows.
\end{proof} 

It remains now to evaluate the probability that the remaining sum in \eqref{lower-slicing}  
be larger than $\varepsilon n$, which we will do in the next subsections. But let us already introduce the notation 
\be{def-xi}
\xi_n(T) = \frac 1T \sum_{k=1}^n \sum_{z\in \kR_k^{++}}  G_T(z-S_k).
\ee
 
\subsection{Case of dimension $d=3$}
The desired estimate on $\xi_n(T)$ will be derived from the following lemma. 

\bl{lem-d3}
There exist positive constants $\eta$, $c$ and $C$, such that for any $\varepsilon\in (0,\nu_3)$, $n\ge 2$, and any $r\le T$, satisfying 
\be{scale-d3}
C\, \varepsilon^{- 5/3}\, n^{2/3}\, \log n\ \le\ r^3 \ \le\  c\, (\varepsilon^{1/2}T^{3/2}\wedge \varepsilon^{3/2} n), 
\ee
the two following statements hold. First
\be{main-ineqd3}
\{\xi_n(T)\ge \varepsilon n \}\ \subseteq \ 
\bigcup_{j=0}^J \bigcup_{x\in \Z^3}
\left\{\left|\kK_n\Big(B(x,r),\eta\varepsilon 2^j \, r^3\Big)\right|\ge 
\frac{n}{2^{\frac{3}{2}j}}\right\},
\ee
with $J$ the smallest integer such that $2^J\ge 1/(\eta \varepsilon)$.
Furthermore, there is $\kappa_3$ such that for any $0\le j\le J$, 
\be{Knalpha}
\pp\left[ \bigcup_{x\in \Z^3} \left\{\left|\kK_n\Big(B(x,r),\eta\varepsilon 2^j\, r^3\Big)\right|\ge 
\frac{n}{2^{\frac 32j}}\right\}\right] \le C\, \exp\left(- \kappa_3\cdot 2^{j/6}\,  \varepsilon^{2/3} n^{1/3} \right). 
\ee
\el 
\br{rem-choiceT}
The upper bound in \eqref{shrink-d3} 
follows from \eqref{lower-slicing}, \eqref{prop-S}, 
and Lemma \ref{lem-d3}. 
Indeed one can just take $T=\varepsilon^{4/3} n^{2/3}$, 
and then any $r$ satisfying \eqref{scale-d3}, 
for instance $r=c^{1/3} \varepsilon^{5/6} n^{1/3}$, 
with $c$ thereof. 
\er

\begin{proof}
For $j\le J$, let 
$$t_j= 8d^2\,\eta \, \varepsilon 2^j r^3 \quad \text{and}\quad  L_j = n 2^{-3j/2},$$
with $\eta$ to be determined later. 
We first observe that, once $\eta$ is fixed, if $r^3 \le c \varepsilon^{3/2} n$, with $c$ small enough, 
then $L_j/t_j\ge 1$, for all $j\le J$. Moreover, it is not difficult to see that if in addition the lower bound in \eqref{scale-d3} holds with $C$ 
large enough, then \eqref{cond-key} is satisfied for $t= t_j$ and $L=L_j$, for all $j\ge 0$. 
Thus, with these choices of $C$ and $c$, \eqref{Knalpha} follows from \eqref{ineq-key} in Lemma~\ref{lem-key} and a union bound. 
Now, for \eqref{main-ineqd3}, we have to show that
when $\eta$ is small enough,
$$
\bigcap_{0\le j\le J}\bigcap_{x\in \Z^3} \left\{\Big| \kK_n\left(B(x,r),
\eta\varepsilon 2^{j} r^3 \right)\Big| <L_j\right\}\subset
\{\xi_n(T)< \varepsilon n\}.
$$
To see this, we rather work with cubes
$Q\in \kP_{r'}$, with $r'=r/\sqrt d-2$, and one assume that $r' \ge 1$, since this is a consequence of \eqref{scale-d3} when $C$ is large enough.    
The reason to introduce $r'$ is that we want  
\be{cube-ball}
\Lambda\subset [-r',r']^d\quad \Longrightarrow\quad 
\Lambda^{++}\subset B(0,r),
\ee
which holds for this choice of $r'$. Now 
for $j \in \{1,\dots, J\}$, define
$$
\kK_n^*(Q, j):= \left\{k\in\{1,\dots,n\}\ :\  t_{j-1} <|\kR_k^{++}
\cap(S_k+Q)|\le  t_j\right\},
$$
and
$$
\kK_n^*(Q, 0):= \left\{k\in\{1,\dots,n\}\ 
:\ |\kR_k^{++}\cap(S_k+Q)|\le t_0\right\}.
$$  
Note that 
$$
\ell_k(Q^{++})\ \ge\ \frac{|\kR_k^{++}\cap Q|}{4d^2},
$$
and combined with \eqref{cube-ball}, this shows that  
$$
\kK_n^*(Q,j)\ \subset\ 
\kK_n\left(B(x,r),\eta\varepsilon 2^{j} r^3\right) 
\qquad \text{for all } j\ge 0.
$$
Note now that by definition of $J$, 
one has $t_J \ge 8r^3\ge |Q|$, for $Q\in \kP_{r'}$ and therefore  
\be{1nKn*}
\{1,\dots,n\}\ \subset \ \bigcup_{j=0}^J\,  \kK_n^*(Q,j).
\ee 
Then, let 
$$
\Xi_n(r,\eta,\varepsilon):=\bigcap_{0\le j\le J}\bigcap_{Q\in \kP_r}
\left\{ | \kK_n^*\left(Q,\ j\right) | <L_j\right\}. 
$$
For $j\in \{0,\dots,J\}$ and $Q\in \kP_{r'}$,  
let $k^*\in \kK_n^*(Q,j)$ be the index 
which maximizes the sum $\sum_{z\in (-S_k+ \kR_k^{++})\cap Q} G_T(z)$,
and define 
$$
\Lambda_Q(j):= (-S_{k^*}+\kR_{k^*}^{++})\cap Q.
$$ 
Note that by definition of 
$\kK_n^*(Q,j)$, one has $|\Lambda_Q(j)|\le \gamma_j r^3$, 
with $\gamma_j = 8d^2 \eta \varepsilon 2^j$. Hence, by using Lemma~\ref{lem-green} and \eqref{1nKn*}, 
we see that on the event $\Xi_n(r,\eta,\varepsilon)$,  
for some constant $A>0$, using the trivial bound $|\kK_n^*(Q,0)|\le n$. 
\be{d3-step1}
\begin{split}
\xi_n(T) \le & \ \frac 1T \sum_{j= 0}^J\sum_{Q\in \kP_{r'}}
\sum_{k\in \kK_n^*(Q,j)}
\sum_{z\in (-S_k+\kR_k^{++})\cap Q} G_T(z)\\
\le & \   \frac 1T\sum_{j= 0}^J\sum_{Q\in \kP_{r'}}|\kK_n^*(Q,j)| \sum_{z\in \Lambda_Q(j)} G_T(z)\\
\le & \ \frac A T\left\{n (\gamma_0^{2/3}r^2
+\gamma_0 T)+\sum_{j= 1}^JL_j \left(\gamma_j^{2/3}r^2 
+ \gamma_j T\right)\right\}.
\end{split}
\ee
Now, if $r^3\le c \varepsilon^{1/2}T^{3/2}$, with $c$ small enough, 
then  
$\gamma_j^{2/3}r^2\le \gamma_j T$, for all $j\ge 0$. Therefore \eqref{d3-step1} shows that if $\eta$ is small enough, then 
\be{xintilde}
\Xi_n(r,\eta,\varepsilon) \ \subset \  \{\xi_n(T)<\varepsilon n\}.
\ee 
Together with \eqref{xintilde} this concludes the proof of the lemma. 
\end{proof}

\subsection{Case of dimensions $d\ge 5$}
The following is an analogue of Lemma \ref{lem-d3}: 
\bl{lem-d5}
Assume that $d\ge 5$. There are positive constants $\eta$, $\kappa_d$, $c$ and $C$, such that 
for any $\varepsilon\in (0,\nu_d)$, $n\ge 2$, and any $r\le T$, satisfying 
\be{scale-d5}
C \varepsilon^{-\frac{d^2-4}{3d}}\, n^{2/d}\, \log n \ \le \ r^d \  \le \  c\, \left(\varepsilon^{\frac{(d -2)^2}6} T^{d/2} \wedge \varepsilon n\right),
\ee
the two following statements hold. First 
\be{main-ineqd5}
\{\xi_n(T)\ge \varepsilon n \}\ \subseteq \ 
\bigcup_{j=0}^J \bigcup_{x\in \Z^d}
\left\{|\kK_n\big(B(x,r), \eta 2^{-j} r^d\big)|\ge  \varepsilon 2^{\frac {2,5j}{d-2}} n \right\},
\ee
with $J$ the smallest integer such that $2^{2,5J/(d-2)}\ge 1/\varepsilon$.
Furthermore, for any $0\le j\le J$, 
\be{Knalpha5}
\pp\left[ \bigcup_{x\in \Z^d} \left\{|\kK_n\big(B(x,r), \eta 2^{-j} r^d\big)|\ge 
\varepsilon 2^{\frac{2,5j}{d-2}}n\right\}\right] \le C\, \exp\left(- \kappa_d\cdot 2^{\frac j{2d}}\,  (\varepsilon n)^{1-\frac 2d} \right).
\ee
\el 
\begin{proof}
Since the proof is entirely similar to the case of dimension $3$, we will not reproduce it here. One just has to use this time  
$t_j= (8d^2) \eta 2^{-j} r^d $ and $L_j =  \varepsilon 2^{2,5j/(d-2)} n$, for $j\le J$. 
\end{proof}
\noindent Note that here, for proving the upper bound in \eqref{shrink-d5}, one can choose $T=\varepsilon^{1+2/d} n^{2/d}$.

\subsection{Case of dimension $d=4$}
In this case we obtain a weaker statement. 
\bl{lem-d4}
Assume that $d=4$. There are positive constants $\eta$, $\kappa_4$, $c$ and $C$, such that 
for any $\varepsilon\in (0,\nu_4)$, $n\ge 2$, and any $r\le T$, satisfying 
\be{scale-d4}
C\,  \varepsilon^{-3/2} |\log \varepsilon |^{3/2}\, n^{1/2}\, \log n \ \le \ r^4\ \le \  c\, \frac{\varepsilon}{|\log \varepsilon|}(T^2\wedge n),
\ee
the two following statements hold. First 
\be{main-ineqd4}
\{\xi_n(T)\ge \varepsilon n \}\ \subseteq \ 
\bigcup_{j=0}^J \bigcup_{x\in \Z^4}
\left\{|\kK_n\big(B(x,r), \eta 2^{-j} r^4\big)|\ge \frac{2^j}J \varepsilon n \right\},
\ee
with $J$ the smallest integer such that $2^J\ge J/\varepsilon$.
Furthermore, for any $0\le j\le J$, 
\be{Knalpha4}
\pp\left[ \bigcup_{x\in \Z^4} \left\{|\kK_n\big(B(x,r), \eta 2^{-j} r^4\big)|\ge 
\frac{2^j}J  \varepsilon n\right\}\right] \le C\, \exp\left(- \kappa_4\cdot \frac{(\varepsilon n)^{1/2}}{|\log \varepsilon|^{1/2}}  \right).
\ee
\el 
\begin{proof}
Since the proof is entirely similar to the previous cases, we leave the details to the reader. Just for \eqref{Knalpha4} one can notice that 
$J$ is of order $|\log \varepsilon |$ by definition. 
\end{proof}

%
\section{Lower bounds}\label{sec-low}
\subsection{Lower bound in dimension $d\ge 4$}
\label{lowerd4}
This is the easiest case. 
We impose that
the walk stays a time $\alpha \varepsilon n$ in a ball of volume $ \varepsilon n/2$, with $\alpha=2/\nu_d$. 
During this short time, the 
boundary of the range is necessarily smaller than $\varepsilon n/2$, and
the rest of the trajectory cannot make up for this loss. 
To be more precise, first write
$$
|\partial\kR_n|\le 
|\partial\kR_{\alpha \varepsilon n}|+|\partial \kR(\alpha \varepsilon n,n)|,
$$
which gives after centering and using \eqref{mean-transient}, 
$$
\overline{|\partial\kR_n|}\le \overline{|\partial\kR_{\alpha \varepsilon n}|}+\overline{
|\partial \kR(\alpha \varepsilon n,n)|}+\kO(\log n).
$$
Let $\rho_n$ be the maximal radius such that $|B(0,\rho_n)|\le \varepsilon n/ 2$. On the event
$\{\kR_{\alpha \varepsilon n}\subset B(0,\rho_n)\}$, one has 
$$|\partial \kR_{\alpha \varepsilon n}|\le |B(0,\rho_n)|\le   \frac{\varepsilon n}{2},$$ 
and therefore also (remembering that $\alpha = 2/\nu_d$)   
\begin{eqnarray*}
\overline{|\partial \kR_n|} \le  -2\varepsilon n+ \frac{\varepsilon n}{2}+\overline{
|\partial \kR(\alpha \varepsilon n,n)|}+\kO(\log n).
\end{eqnarray*}
Hence, for $n$ large enough
\be{LB-4}
\pp\big[\overline{|\partial\kR_n|}\le -\varepsilon n\big]\ge 
\pp\big[\kR_{\alpha \varepsilon n}\subset B(0,\rho_n)\big]-
\pp\Big[\overline{|\partial \kR(\alpha \varepsilon n,n)|}\ge   \frac 14 \varepsilon n\Big].
\ee
Then one can use Okada's results \cite{Ok1} which show that upper large deviations have exponentially small probability. 
More precisely, his results imply in particular that 
for any $\varepsilon>0$, there exists a constant $c=c(\varepsilon)>0$, such that 
\be{ULD}
\pp\left[\overline{|\partial \kR_n|}\ge \varepsilon n\right] \le \exp(-c n),
\ee
for $n$ large enough.  On the other hand it is well known that there exists a constant $\kappa>0$, such that for any $r\ge 1$ and $n\ge 1$, 
\be{DV-inf}
\pp\big[\kR_n\subset B(0,r)\big] \ge \exp\left(-\kappa\cdot \frac{n}{r^2}\right).
\ee
Combining \eqref{LB-4}, \eqref{ULD} and \eqref{DV-inf} gives the lower bounds in \eqref{shrink-d5} and \eqref{shrink-d4}. 

\subsection{Lower bound in dimension $d=3$}
This case is more delicate than the previous one. 
We distinguish two regimes: 
when $\varepsilon$ is small, say $\varepsilon \le \varepsilon_0$, 
for some small enough constant $\varepsilon_0$, 
and when $\varepsilon \in (\varepsilon_0,\nu_3/2)$. 
The latter case can be handled using the same argument 
as in higher dimensions, and we omit it, and rather
concentrate on the former.  

We will see that to shrink the boundary of the range, 
the random walk needs to localize a time $n$ 
in a ball of radius $\rho_n$ with $\rho_n^3$ of order $n/\varepsilon$. 
The heuristics behind this picture relies on the following key relation already mentioned in Subsection \ref{subsec.RBR}: 
for any integers $n$ and $m$
\be{main-amin}
|\partial \kR(0,n+m)|\le |\partial \kR(0,n)|+|\partial \kR(n,n+m)| -
|\kR(0,n)\cap \partial \kR(n,n+m)|.
\ee
Without constraint, the intersection of two
strands $|\kR(0,n)\cap \partial \kR(n,2n)|$ is typically of order
$\sqrt n$ in dimension three, 
and does not influence the (linear) growth of the boundary
of the range. However, when the walk is localized in a ball of volume
$n/\varepsilon$, it likely visits an $\varepsilon$-fraction
of all fixed {\it large} volume -- see Proposition \ref{prop-cover} below for a precise statement -- so that the former intersection is
typically of order $\varepsilon n$ when we choose $m=n$, and realize
that $|\partial \kR(n,2n)|$ is typically of order $n$. 
The fact that this scenario leads to the correct cost follows from 
the bound \eqref{DV-inf}, which was already instrumental 
in higher dimensions. 

\bigskip 
\noindent Let us give now some details. The proof is based on the following covering result.
\bp{prop-cover}
There are positive constants $c$, $C$, and $\varepsilon_0$, 
such that for all $\varepsilon \in (0,\varepsilon_0)$, and
$n$ large enough, we have
\be{need-5}
\pp\Big[|\kR_n\cap \Lambda|> \varepsilon |\Lambda |\Big]\, \ge\,  \exp(-\,  \varepsilon^{2/3}n^{1/3})\qquad \textrm{for all } \Lambda\, \subset\,  B\Big(0,c(\frac{n}{\varepsilon})^{1/3}\Big) \textrm{ with }
|\Lambda| \ge \frac{C}{\varepsilon^3}.
\ee
\ep
\noindent Before proving this result, 
let us deduce the desired lower bound. 

\vspace{0.2cm}
\noindent{\bf Proof of the lower bound in \eqref{shrink-d3}:}
We center the variables of \eqref{main-amin} using
\eqref{mean-transient} to obtain 
$$
\overline{|\partial\kR(0,2n)|}\le \overline{|\partial
\kR(0,n)|}+\overline{|\partial \kR(n,2n)|}-
|\kR(0,n)\cap\partial\kR(n,2n)|+ \kO(\sqrt n).
$$
By invariance of time-inversion and using the Markov property, 
we see that the intersection term in this inequality is equal in law to the intersection of two independent ranges. 
Hence Okada's large deviations estimate \eqref{ULD} yields,
for $n$ large enough, 
\begin{eqnarray*} 
\pp[\overline{|\partial\kR_{2n}|}\le - 2 \varepsilon n]&\ge & \pp[|\kR_n\cap\partial\widetilde \kR_n|\ge 4\varepsilon n] - 
2\pp[\overline{|\partial \kR_n|}\ge \varepsilon n/2]\\
&\ge & \pp[|\kR_n\cap\partial \widetilde \kR_n|\ge 4\varepsilon n]- 
2\exp(-c n),
\end{eqnarray*}
with $\widetilde \kR_n$ an independent copy of $\kR_n$, 
and $c=c(\varepsilon)>0$ a constant.
We are therefore left with showing the existence of $\kappa>0$ 
such that for all $\varepsilon$ small enough and $n$ large enough 
\be{need-3}
\pp\big[|\kR_n  \cap \partial \widetilde \kR_n|\ge \varepsilon n]\ge
\exp(-\kappa\,  \varepsilon^{2/3} n^{1/3}).
\ee
To this end, we first claim that localizing a random walk a time $n$ inside $B(0,\rho_n)$ with $\rho_n^3$ of order $n/\varepsilon$ 
still produces a boundary of the range of order $n$. Indeed,
$$
\pp\big[|\partial \widetilde \kR_n|\ge \nu_3 n/2,\ 
\widetilde \kR_n\subset B(0,\rho_n)\big]\ge
\pp\big[\widetilde \kR_n\subset B(0,\rho_n)\big]-\pp\big[|\partial \widetilde  \kR_n|\le \nu_3 n/2\big].
$$
Then using \eqref{DV-inf}, we get  
\be{DV-inf.bis}
\pp\big[\widetilde \kR_n\subset B(0,\rho_n)\big]\ge \exp(-\kappa\, \varepsilon^{2/3} n^{1/3}),
\ee
for some constant $\kappa>0$ (here we take $\rho_n=c(n/\varepsilon)^{1/3}$, with $c$ as in Proposition \ref{prop-cover}) and from our upper bound, 
we have for some constant $\kappa'>0$, 
\[
\pp\Big[|\partial \widetilde \kR_n|\le \nu_3 n/2\Big]\le \exp(- \kappa'\, n^{1/3}).
\]
Thus, if we take $\varepsilon$ small enough, we see that  
for some possibly larger constant $\kappa>0$
\be{loc+boundary}
\pp\Big[|\partial \widetilde \kR_n|\ge \nu_3 n/2,\ \widetilde \kR_n\subset B(0,\rho_n)\Big]\ge \exp(-\kappa \, \varepsilon^{2/3} n^{1/3}).
\ee
Then by using the independence between $\kR_n$ and $\widetilde \kR_n$, 
Proposition~\ref{prop-cover}, and \eqref{loc+boundary} we have
\begin{eqnarray*}
\pp\Big[|\kR_n\cap\partial \widetilde \kR_n|\ge  \varepsilon n \Big] 
&\ge& 
\sum_{V\subset B(0,\rho_n),\, |V|\ge \nu_3 n/2} \pp \big[|\kR_n\cap V|\ge  \varepsilon n\big] \times \pp[\partial \widetilde \kR_n = V] \\
& \ge &  \exp(-\, \varepsilon^{2/3} n^{1/3})\times \pp\big[|\partial \widetilde \kR_n|\ge \nu_3 n/2,\ \widetilde \kR_n\subset B(0,\rho_n)\big]\\
&\ge & \exp(-(1+\kappa)\, \varepsilon^{2/3} n^{1/3}).
\end{eqnarray*}
This proves \eqref{need-3}. 
\qed

\vspace{0.2cm}
\noindent We are left with proving Proposition~\ref{prop-cover}.

\vspace{0.2cm}
\noindent{\bf Proof of Proposition~\ref{prop-cover}.} 
Define $\rho_n$ and $K_n$ by $\rho_n=c(n/\varepsilon)^{1/3}$, and  
$K_n=K \varepsilon \rho_n$, with 
$c$ and $K$ some constants to be chosen later.
Observe that a walk covers a fraction $\varepsilon$ of 
a set $\Lambda\subset B(0,\rho_n)$, 
if it makes $K_n$ excursions between $B(0,2\rho_n)$
and $\partial B(0,5\rho_n)$ before time $n$.
Indeed, each excursion has a chance of order $1/\rho_n$ 
to visit any given site of $B(0,\rho_n)$ (since $\rho_n$ is the typical distance between such site and $\partial B(0,2\rho_n)$), 
so that $K_n$ 
{\it independent} excursions have a chance of order $K\varepsilon$
to cover any given site. 
If $K$ is large enough, one deduces well that $K_n$ excursions cover at least a fraction $\varepsilon$ of $\Lambda$. 
Now, a Gambler's ruin estimate shows that the probability 
to hit the ball $B(0,2\rho_n)$ before exiting $B(0,8\rho_n)$ starting from 
$\partial B(0,5\rho_n)$ is bounded away from $0$, so
this happens $K_n$ times at a cost 
$\exp(-K_n)$ which is of the right order.

Finally note that the length of a typical excursion is of order $\rho_n^2$, 
so that with high probability $K_n$ of them
have to occur before time $n$, provided $c$ is small enough, and this 
concludes the heuristics of the proof. 

\vspace{0.2cm}
Let us proceed now to the details. We first define the excursions in a standard way as follows. 
Set $\sigma_0=0$, and then for $i\ge 0$, 
$$\tau_i:=\inf\{t\ge \sigma_i\ :\ S_t\notin B(0,5 \rho_n)\},$$
and 
$$\sigma_{i+1}:=\inf\{t\ge \tau_i\ :\ S_t\in B(0,2\rho_n)\}.$$
Also define $\tau^*$ to be the exit time from $B(0,8\rho_n)$.

The number of excursions (after $\tau_0$) from $\partial B(0,2\rho_n)$ to $\partial B(0,5\rho_n)$ is defined to be 
$$
\kN:= \sup\ \{k\ge 0\ :\ \sigma_k<\infty\}.
$$
Recall that $K_n=K\varepsilon \rho_n$, for
$K$ to be fixed later. Let $\kG_{\kN}$
be the sigma-field generated by $\kN$,
and the starting and end points of the excursions.  
$$
\kG_{\kN}:=
\sigma\left(\kN,\, X_{\sigma_i},X_{\tau_i},\ i\le 
\kN\right).
$$
Define $\Lambda_1=\Lambda$, and for $i\ge 1$, define 
$$
\kR^{(i)} = \{X_{\sigma_i},\dots,X_{\tau_{i}}\}
\quad \text{and}\quad \Lambda_{i+1}= 
\Lambda\setminus (\cup_{j\le i} \kR^{(j)}\cap \Lambda).$$
Then, set
\be{def-X}
X_i:= |\kR^{(i)}\cap \Lambda_i|\, \ind_{\{\sigma_i<\infty\}}.
\ee 
Notice that conditionally on $\kG_{\kN}$, 
the variables $\{X_i,\ i\le \kN\}$ are independent 
of the event $\{\sigma_{K_n}\le \tau^*\}$. Therefore,
\be{egalite.esp.GN}
\pp\left[\sum_{i=1}^{K_n} X_i >\varepsilon |\Lambda|, \, 
\sigma_{K_n} \le \tau^*\right] 
=\ \E\left[ \ind_{\{\kN\ge K_n\}}
\pp\Big[\sum_{i=1}^{K_n}X_i>\varepsilon |\Lambda|\ \Big| \  \kG_\kN\Big] 
\pp\Big[\sigma_{K_n} \le \tau^*\mid \kG_\kN\Big] \right].
\ee
Let $\kH_i$ the sigma-field generated by the walk up to 
the stopping time $\sigma_i$. Define  
\be{Mn.N}
M_n =
\sum_{i=1}^{\kN \wedge K_n} (X_i - \E[X_i\mid \kH_i,\, \kG_\kN]),
\ee
and note that
\be{Mn.N2}
\E[M_n\ \big|\ \kG_\kN]=0 \quad\text{and}\quad 
\E[M_n^2\ \big|\ \kG_\kN]\le 2\sum_{i=1}^{\kN\wedge K_n} 
\E[X_i^2\mid\ \kG_\kN].
\ee
Now, for any $i\le \kN$, we have 
\be{esp.Xi.cond}
\E[X_i\mid \kH_i,\, \kG_\kN]= \sum_{x\in \Lambda_i} \pp[x\in \kR^{(i)}\mid \kG_\kN]=\sum_{x\in \Lambda_i} \pp\left[x\in \kR^{(i)}\mid X_{\sigma_i},\, X_{\tau_i}\right],
\ee
since (conditionally on $\kG_\kN$), the law of $X_i$ depends on $\kH_i$ only through $\Lambda_i$. 
Moreover, as a consequence of the Harnack principle 
there exists a constant $c_0>0$, such that 
for any $x\in B(0,\rho_n)$,
\be{claim.proba.cond}
\pp[x\in \kR^{(i)}\mid X_{\sigma_i},\, 
X_{\tau_i}]\ge c_0\, \pp[x\in \kR^{(i)}\mid X_{\sigma_i}]. 
\ee
Indeed, for any $y\in \partial B(0,5\rho_n)$, the Markov property yields
$$
\pp[x\in \kR^{(i)}\mid X_{\sigma_i},\, X_{\tau_i}=y]  =  \pp[x\in \kR^{(i)}\mid X_{\sigma_i}] \times \pp_x\left[X_{H_{\partial B(0,5\rho_n)}}=y\right],
$$
and now the Harnack principle, see \cite[Theorem 6.3.9]{LL}, shows that there exists $c_0>0$, independent of $x$, $y$ and $X_{\sigma_i}$, such that almost surely, 
$$\pp_x[X_{H_{\partial B(0,5\rho_n)}}=y] \ge c_0\,  \pp[X_{\tau_i} = y \mid X_{\sigma_i}],$$
which proves \eqref{claim.proba.cond}. Then for any $z \in \partial B(0,2\rho_n)$, by using standard estimates on the Green's function (see \cite[Theorem 4.3.1]{LL}), we can write   
\begin{equation*}
\begin{split}
\pp[x\in \kR^{(i)}\mid X_{\sigma_i}=z] & \ge \  \pp_z[H_x<\infty] - \sup_{y\in \partial B(0,5\rho_n)} \pp_y[H_x<\infty] \\
& \ge\  \frac {c_{\text{Gr}} }{\|x-z\|} - \sup_{y\in \partial B(0,5\rho_n)} \frac {c_{\text{Gr}} }{\|y-x\|} - \kO\left(\frac 1{\rho_n^2}\right)\\
& \ge \ \frac {c_{\text{Gr}} }{3\rho_n} - \frac {c_{\text{Gr}} }{4\rho_n} - \kO\left(\frac 1{\rho_n^2}\right)\\ 
& \ge \ \frac{c_{\text{Gr}} }{20\rho_n},
\end{split}
\end{equation*}
for some constant $c_{\text{Gr}}>0$, and $n$ large enough. Combining this with \eqref{esp.Xi.cond} and \eqref{claim.proba.cond}, we obtain that for any $i\le \kN$, almost surely 
\be{esp.Xi.cond.2}
\E[X_i\mid \kH_i,\, \kG_\kN]\ge \frac {c_1}{\rho_n}|\Lambda_i|, 
\ee
for some constant $c_1>0$. Now, choose $K:=4/c_1$, and use 
the previous inequality. This shows that on the event $\{\kN\ge K_n\}$,
$$
\sum_{i=1}^\kN \E[X_i\mid \kH_i,\, \kG_\kN]
\ge 4\varepsilon |\Lambda_{K_n}|.
$$
Since $|\Lambda_{K_n}|=|\Lambda|-\sum_{i=1}^{K_n-1} X_i$,
we have for any $\varepsilon \le 1/2$, on the event $\{\kN\ge K_n\}$,
\be{mart.ineq}
\begin{split}
\pp\left[\sum_{i=1}^{K_n} X_i \le \varepsilon |\Lambda| \ \Big| \  \kG_\kN\right]&  =
\ \pp\left[\sum_{i=1}^{K_n} X_i \le \varepsilon |\Lambda|,\, |\Lambda_{K_n}|
\ge |\Lambda|/2 \ \Big| \  \kG_\kN\right]\\
&\le \ \pp[|M_n|\ge \varepsilon |\Lambda| \mid \kG_\kN]\\
&\le \frac{\E\left[M_n^2\mid \kG_\kN\right]}{\varepsilon^2 |\Lambda|^2}\\
&\le 
\frac{2}{\varepsilon^2 |\Lambda|^2} 
 \sum_{i=1}^{K_n} 
\E\left[X_i^2\mid \kG_\kN\right],
\end{split}
\ee
using \eqref{Mn.N2} at the last line. 
Moreover, for any $i\le \kN$, using again the Harnack inequality at the third line, and \eqref{Green} at the last line, 
\be{carre.Xi}
\begin{split}
\E\left[X_i^2\mid \kH_i, \, \kG_\kN\right] &=\  \sum_{(z,z')\in \Lambda_i\times \Lambda_i}\pp\left[z\in \kR^{(i)},\, z'\in \kR^{(i)}\mid X_{\sigma_i},\, X_{\tau_i}\right] \\
& \le \ 2 \, \sum_{(z,z')\in \Lambda_i\times \Lambda_i}\pp\left[z\in \kR^{(i)},\, z'\in \kR^{(i)},\, H_z<H_{z'}\mid X_{\sigma_i},\, X_{\tau_i}\right] \\
&\le \ C_0\, \sum_{(z,z')\in \Lambda_i\times \Lambda_i}\pp\left[z\in \kR^{(i)},\, z'\in \kR^{(i)},\, H_z<H_{z'}\mid X_{\sigma_i}\right] \\ 
&\le \ C_0\, \sum_{(z,z')\in \Lambda_i\times \Lambda_i}\pp_{X_{\sigma_i}}\left[H_z<H_{z'}<\infty\right]
\le \ C_0 \sum_{z,z'\in \Lambda_i} \frac{1}{\rho_n\, (\|z-z'\|+1)},
\end{split}
\ee 
for some constant $C_0>0$. 
Now it is not difficult to see that for any set $\Lambda$, and any $y\in \Z^d$, 
one has $\sum_{z\in \Lambda} 1/(\|z-y\|+1)
 =\kO(|\Lambda |^{2/3})$, 
with an implicit constant which is uniform in $\Lambda$ and $y$, 
since the worst case is easily seen to be reached 
when $\Lambda$ is a ball and $y$ is the center of this ball. Therefore \eqref{carre.Xi} gives 
$$\E\left[X_i^2\mid \kH_i, \, \kG_\kN\right] \le \frac{C'_0}{\rho_n}\, |\Lambda|^{5/3},$$  
for some constant $C'_0>0$. 
Together with \eqref{mart.ineq} 
if we assume that $|\Lambda |\ge C/\varepsilon^3$, 
with $C:=(4KC'_0)^3$, one deduces that, on the event $\{\kN\ge K_n\}$, 
$$ 
\pp\left[\sum_{i=1}^{K_n} X_i \le 
\varepsilon |\Lambda| \ \Big| \  \kG_\kN\right]
\le \frac{2KC'_0}{\varepsilon |\Lambda|^{1/3}}\le \frac 12.$$
Coming back to \eqref{egalite.esp.GN} we get 
\be{inegalite.esp.GN}
\begin{split}
\pp\left[\sum_{i=1}^{K_n} 
X_i >\varepsilon |\Lambda|, \, \sigma_{K_n} \le \tau^*\right] 
&\ge \  \frac12 \, \pp\Big[\sigma_{K_n}\le \tau^*\Big] \\
&\ge \ \frac12 \, \exp(-\kappa\cdot K c^{1/3} \cdot \varepsilon^{2/3}n^{1/3} ), 
\end{split}
\ee
for some constant $\kappa>0$, since for any $x\in \partial B(0,5\rho_n)$ the probability starting from $x$ to hit the ball $B(0,2\rho_n)$ before $\tau^*$ is bounded away from $0$ (see \cite[Proposition 1.5.10]{L}). 
Now, observe that
\begin{equation}
\label{last.step.min} 
\begin{split}
\pp[|\kR_n\cap \Lambda|\ge \varepsilon |\Lambda|] & \ge \ \pp\left[
\sum_{i=1}^{K_n} X_i >\varepsilon |\Lambda|, \, \sigma_{K_n}<\tau^*,
\, \tau_{K_n} \le \tau^*\wedge n\right]\\  
&\ge \ \pp\left[\sum_{i=1}^{K_n} X_i >\varepsilon |\Lambda|, 
\, \sigma_{K_n} \le \tau^*\right] - \pp[n<\tau_{K_n} \le \tau^*].
\end{split}
\end{equation}
We now bound the last term in the right hand side of \reff{last.step.min}.
To this end we let 
$$\xi_i = (\sigma_{i+1}\wedge \tau^* - \tau_i)\ind_{\{\tau_i<\infty\}}  \quad \text{and}\quad \tilde \xi_i =
(\tau_i - \sigma_i)\ind_{\{\sigma_i<\infty\}}\quad \text{for }i\ge 0,$$
and observe that for any integer $k\ge 1$, on the event $\{\tau_k\le \tau^*\}$, we have  
\be{sigmakxii}
\tau_k =\tilde \xi_k+ \sum_{i=0}^{k-1} (\xi_i + \tilde \xi_i).
\ee
Then we claim that the random variables $(\xi_i/\rho_n^2)$ and $(\tilde \xi_i/\rho_n^2)$ are (up to a multiplicative constant) dominated by i.i.d. geometric random variables. 
To see this, one can use that $(\|S_n\|^2-n)_{n\ge 0}$ is a martingale. It implies, using also the optional stopping time theorem, 
that for any $x\in B(0,8\rho_n)$, 
$$(8\rho_n+1)^2\ge \E_x\left[\|S_{\tau^*\wedge 200\rho_n^2}\|\right] \ge \E_x\left[\tau^*\wedge 200\rho_n^2\right] \ge 200\rho_n^2\,  \pp_x[\tau^*\ge 200\rho_n^2].$$
Hence for any $x\in B(0,8\rho_n)$, 
$$\pp_x[\tau^*<200\rho_n^2]\ge 1/2.$$
We deduce that there exist $(G_i)$ and $(\tilde G_i)$, i.i.d. geometric random variables with parameter $1/2$, such that 
$$\xi_i/(200\rho_n^2) \le G_i \quad \text{and}\quad \tilde \xi_i/(200\rho_n^2)\le \tilde G_i\quad \text{for all }i\ge 0.$$
Now assume that $c<1/(2000K)$, so that $1/(200Kc)\ge 10$. 
Then with \eqref{sigmakxii}, and using Markov's exponential inequality, this gives 
\begin{equation*}
\begin{split}
\pp[n< \tau_{K_n}\le \tau^*] & \le \  \pp\left[\sum_{i=0}^{K_n} (G_i+\tilde G_i)\ge  \frac{n}{200\rho_n^2} \right]\\
& \le \  \pp\left[\frac{1}{K_n} \sum_{i=0}^{K_n} (G_i+\tilde G_i)\ge  \frac{1}{200Kc} \right]\\
&\le \ \exp(-\kappa\cdot c^{-2/3} \cdot \varepsilon^{2/3} n^{1/3}), 
\end{split}
\end{equation*}
for some constant $\kappa>0$. 
Finally by taking $c$ small enough and using \eqref{inegalite.esp.GN} and \eqref{last.step.min} this proves the desired lower bound. 
\qed


\section{Bound on the capacity}\label{sec-cap}
In this Section, we prove Proposition \ref{prop.cap}, and then 
use it for proving \eqref{new-d3}, \eqref{new-d5} 
and \eqref{new-d4}.

We first recall some important property of the capacity of a ball (see (2.16) in \cite{L} for a stronger statement): there is a positive constant 
$C_{\capa}$, such that 
for all $x\in \Z^d$ and $r>0$ 
\be{cap.boules}
\capa(B(x,r)) \, \le \, C_{\capa} \ |B(x,r)|^{1-2/d}.
\ee

\paragraph{Proof of Proposition~\ref{prop.cap}.}
First, consider the case when $\kI_d(B(\kC,r))$ is bounded above by some constant $\overline C$. 
Then by using \reff{Aminebound} (see \cite[Lemma 1.2]{AC}),
we obtain 
\be{old-ineq}
\begin{split}
\pp\big[ \ell_n(B(x,r))\ge t\ \textrm{ for all }x\in \kC \big]\ \le &\  
\pp\big[\ell_n(B(\kC,r))\ge t|\kC|\big]\\
\le  &\ C\, \exp\Big(-\kappa_0\cdot  \frac{|\kC|t}{|B(\kC,r)|^{2/d}}\Big),
\end{split}
\ee
which proves the desired inequality \reff{main-tech}, with $\kappa = \kappa_0/\overline C$. Henceforth, we assume that
\be{old-cond}
\frac{\capa(B(\kC,r))}{|\kC|^{1-2/d}\cdot r^{d-2}}> 2 C_{\capa}.
\ee
We need now a few intermediate lemmas. 
We first show that spending a time $t$ in a ball $B(x,r)$ implies making order $t/r^2$ excursions from $\partial B(x,r)$ to $\partial B(x,2r)$, with high probability. So as in the proof of Proposition~\ref{prop-cover}, set $\sigma_0(x,r)=0$, and for $j\ge 0$, 
$$\tau_j(x,r):=\inf\{t\ge \sigma_j(x,r)\ :\ S_t\in \partial B(x,2r)\},$$
and 
$$\sigma_{j+1}(x,r):= \inf\{t\ge \tau_j(x,r)\ :\ S_t\in B(x,r)\}.$$
Then for $n\ge 1$, let 
$$N_n(x,r):=\sup \{j \ :\ \sigma_j(x,r)\le n\}.$$
\bl{lem.exc}
There exists a constant $c>0$, such that for any $r\ge 1$, 
$\kC\in \kA(r)$, and $n\ge 1$,
$$
\pp\left[
\ell_n(B(x,r))\ge t\text{ and } N_n(x,r) \le c t/r^2\ \textrm{for all }x\in \kC\right]
\le \exp\big(-c\, |\kC|\, \frac{t}{ r^2}\big).
$$
\el
\begin{proof}
Exactly as in the proof of Proposition~\ref{prop-cover}, 
one can see that there exist i.i.d. geometric random variables 
$\{G_j(x),\ j\in \N, x\in \kC\}$, with parameter $1/2$, 
such that for any $j\in \N$, and $x\in \kC$,  
$$(\tau_j(x,r)-\sigma_j(x,r))\ind_{\{\sigma_j(x,r)<\infty\}} \le 100r^2 G_j(x).$$ 
Moreover, for any $x\in \kC$, 
$$\ell_n(B(x,r))\le \sum_{j=0}^{N_n(x,r)} (\tau_j(x,r)-\sigma_j(x,r)).$$
Indeed, we include the 0-th order term $\tau_0(x,r)$ to cover
the case where the walk starts in $B(x,2r)$.
Therefore, by standard large deviation estimates, 
for any $\delta\in (0,1/400)$, there is $\gamma>0$ (independent
of $n,t$ and $r$), such that
\begin{equation*}
\begin{split}
\pp\left[
\ell_n(B(x,r))\ge t,\ N_n(x,r) \le \delta\frac{t}{r^2}
\textrm{ for all }x\in \kC\right] &\le
\pp\left[\sum_{j=0}^{\delta t / r^2} G_j(x) 
\ge \frac{t}{100r^2} \textrm{ for all }x\in \kC \right]\\
&\le \prod_{x\in \kC}
\pp\left[\sum_{j=0}^{\delta t / r^2} G_j(x)
\ge \frac{1}{200\delta} \frac{2\delta t}{r^2}\right] \\
&\le  \exp(- \gamma |\kC| \frac{t}{r^2}).
\end{split}
\end{equation*}
The result follows as we take $c=\min(\delta,\gamma)$. 
\end{proof} 
\noindent The next result relates the probability to never return 
to $B(\kC,r)$ starting from a boundary point 
of $B(x,2r)$ to the probability of 
the same event starting from a uniformly chosen site on $\partial B(x,r)$.  

\bl{equiv.hitting}
With the notation and hypothesis of Proposition \ref{prop.cap}, 
there exists $\theta>0$ (independent of $r$ and $\kC\in \kA(r)$), 
such that for any $x\in \kC$ and $z\in \partial B(x, 2r)$, 
$$\pp_z[H_{B(\kC,r)}= +\infty]\ge  \theta\cdot \frac{1}{r^{d-2}} \sum_{y\in \partial B(x,r)} \pp_y[H^+_{B(\kC,r)}= +\infty].$$
\el
\begin{proof} We first argue 
that there exists $\theta_1>0$, such that 
for any $z'\in \partial B(x,2r)$
\be{theta1}
\pp_z[H_{B(\kC,r)}= +\infty]\ \ge\  \theta_1 \cdot \pp_{z'}[H_{B(\kC,r)}= +\infty].
\ee
To this end, let $\tau$ be the exit time from the annulus
$B(x,5r/2) \bs B(x,3r/2)$.  
By using the optional stopping time theorem, we obtain
$$\pp_z[H_{B(\kC,r)}= +\infty] = \sum_{v}
\,  \pp_z[S_\tau = v]\, \pp_v[H_{B(\kC,r)} = +\infty].$$
Then, Harnack's inequality (see \cite[Theorem 6.3.9]{LL}) 
shows that there exists a constant $\theta_1$, such that 
for all $z'\in \partial B(x,2r)$ and $v$, 
$$\pp_z[S_\tau = v]\ \ge \ \theta_1\cdot \pp_{z'}[S_\tau = v]. $$
Inequality \eqref{theta1} follows. 
Now using Proposition 1.5.10 in \cite{L}, we see that there exists $\theta_2>0$, such that 
\be{theta2}
\pp_y[H_{\partial B(x,2r)} < H_{\partial B(x,r)}^+] \ \le\  \theta_2 \, r^{-1},
\ee
for all $y\in \partial B(x,r)$. Next, using Lemma 6.3.7 and Proposition 6.4.4 in \cite{LL}, we deduce that there exists $\theta_3>0$ satisfying 
for all $y\in \partial B(x,r)$ and $w\in \partial B(x, 2r)$, 
\be{theta3}
\pp_y[S_{\tau'} = w]\ \le\  \theta_3 \cdot \frac{1}{r^{d}} \qquad \textrm{with } \tau' = H_{ \partial B(x, 2r) \cup B(x,r)}^+. 
\ee
Then by using that the size of the boundary of a ball 
of radius $r$ or $2r$ is of order $r^{d-1}$, we deduce that 
$$\frac{1}{|\partial B(x,2r)|} \sum_{w\in \partial B(x,2r)}\pp_w[H_{B(\kC,r)} = +\infty]\ \ge \ \theta_4 \cdot \frac{1}{r^{d-2}}\sum_{y\in \partial B(x,r)}\pp_y[H^+_{B(\kC,r)} = +\infty].$$ 
Combining this with \eqref{theta1} gives the result.    
\end{proof}

\bl{lem.rec.cap}
For $x\in\kC$, define 
$$
q_x = \exp\left(-\min_{z\in \partial B(x,2r)} 
\pp_z[H_{B(\kC,r)}= +\infty]\right).$$
Then, for any $x_0\in \kC$, any $y\in \partial B(x_0,2r)$, 
and any integers $n$ and $(n_x,x\in \kC)$, 
\be{recurrence.cap}
\pp_y\left[N_n(x,r)= n_x\ \textrm{ for all }x\in \kC \right]\ 
\le \ \frac{q_{x_0}}{\min_{x\in\kC} q_x} \ \prod_{x\in \kC} q_x^{n_x}.
\ee
As a consequence, we have
\be{recurrence.cap2}
\pp_y\left[N_n(x,r) \ge  n_x \ \textrm{ for all }x\in \kC \right]\ \le \ n^{|\kC|}
\, \frac{q_{x_0}}{\min_x q_x} \ \prod_{x\in \kC} q_x^{n_x}.
\ee
\el
\begin{proof}
We only need to prove \eqref{recurrence.cap}, 
since \eqref{recurrence.cap2} follows immediately using that in any ball $B(x,r)$, there can be at most $n$ excursions before time $n$. So we prove \eqref{recurrence.cap} by induction on 
\[
N:=\sum_{x\in \kC} n_x.
\]
First assume that $N=0$. Then we just bound the left-hand side
in \eqref{recurrence.cap} by $1$, and observe that the right-hand side is equal to $q_{x_0}/\min_x q_x$, which is well always larger than or equal to one. Next assume that $N=1$. In this case the left-hand side in \eqref{recurrence.cap} is bounded above by 
$$\pp_y[H_{B(\kC,r)}\le n]\le 1-\pp_y[H_{B(\kC,r)}=+\infty]\le \exp(-\pp_y[H_{B(\kC,r)}=+\infty])\le q_{x_0},$$
which proves the case $N=1$. 
We now prove the induction step, and assume for this that $N\ge 2$. 
We 
write $P_{n,y}(n_x,\ x\in \kC)$ for 
the term on the left-hand side of \eqref{recurrence.cap}, and define 
$$\tau = \inf \{t\ge H_{B(\kC,r)}\ :\ S_t\in \partial B(\kC,2r)\}.$$ 
Now, by the Markov property 
$$
P_{n,y}(n_x,\ x\in \kC)=\sum_{x'\in \kC}
\sum_{k=1}^n \sum_{z\in \partial B(x', 2r)} 
\pp_y\left[\tau=k, \ S_\tau= z\right] \ 
P_{n-k,z}(n_x-\ind_{\{x=x'\}},\ x\in \kC),
$$
where the first sum is over all centers $x'$, 
such that $n_{x'}\ge 1$. By using the induction hypothesis, this gives 
$$
P_{n,y}(n_x,\ x\in \kC)\le 
\frac{\prod_{x\in \kC} q_x^{n_x}}
{\min_{x\in \kC} q_x}\times \pp_y[H_{B(\kC,r)} \le n],$$
and we conclude as in the case $N=1$.   
\end{proof}
\noindent We are now in position to give the proof of Proposition \ref{prop.cap}. 

\vspace{0.2cm}
\noindent \textit{Proof of Proposition \ref{prop.cap}}. 
Let $k_0$ be the integral part of 
$\text{cap}(B(\kC,r)) /(2C_{\text{cap}}r^{d-2})$, with 
$C_{\text{cap}}$ as in \eqref{cap.boules}. Note that by 
\reff{old-cond}, one has $k_0\ge 1$, and therefore 
\be{k0}
\frac{1}{4 C_{\text{cap}}} \cdot 
\frac{\text{cap}(B(\kC,r))}{r^{d-2}}\, \le \, k_0\, \le\,  \frac{1}{2 C_{\text{cap}}} \cdot \frac{\text{cap}(B(\kC,r))}{r^{d-2}}.
\ee
Now, for $x\in \kC$ define
$$\text{cap}_x(B(\kC,r)):= 
\sum_{y\in \partial B(x,r)}\pp_y[H^+_{B(\kC,r)}=+\infty].$$
Using \eqref{cap.boules}, one can see that for any $x\in \kC$, 
$$\text{cap}_x(B(\kC,r)) \le \text{cap}(B(x,r))\le C_{\text{cap}}\ r^{d-2}.$$
Moreover, by definition, 
$$\text{cap}(B(\kC,r))=\sum_{x\in \kC}\text{cap}_x(B(\kC,r)).$$ 
Therefore for any subset $I\subset \kC$ 
with cardinality $m-k_0$, one has 
\be{dernier.prob.cap}
\sum_{x\in I} \text{cap}_x(B(\kC,r)) \ge \text{cap}(B(\kC,r))-
k_0  C_{\text{cap}}\ r^{d-2} \ge \ \frac{\text{cap}(B(\kC,r))}{2}.  
\ee
Next let 
$$E=\left\{\ell_n(B(x,r))\ge t \textrm{ for all }x\in \kC\right\}, $$
and 
$$F=\left\{N_n(x,r)\ge ct/r^2 \textrm{ for at least }|\kC|-k_0 \text{ elements }x\in \kC\right\}, $$ 
with $c$ as in Lemma \ref{lem.exc}. 
By using Lemma \ref{lem.exc} we get for some positive constant $C$, 
\be{prop.cap.1}
\pp\left[E\cap F^c\right] \ \le\  C\cdot {|\kC| \choose k_0} 
\cdot \exp\left(-c\cdot k_0\cdot \frac{t}{r^2}\right).
\ee 
On the other hand Lemma \ref{equiv.hitting}, 
\eqref{recurrence.cap2} and \eqref{dernier.prob.cap} show that 
for some positive constant $C$, 
with $\theta$ as in Lemma \ref{equiv.hitting}, 
\be{prop.cap.2}
\pp[F]\  \le\  C\cdot {|\kC|\choose k_0}\cdot n^{|\kC|-k_0}\, 
\exp\left(- \theta c \cdot \frac{\text{cap}(B(\kC,r))}{2} 
\cdot \frac{t}{r^d}\right).
\ee
Then, \eqref{k0}, \eqref{prop.cap.1} and \eqref{prop.cap.2} 
prove the proposition. \hfill $\square$

\paragraph{Proof of \eqref{new-d3}, \eqref{new-d5} and \eqref{new-d4}.} 
We only give the proof of \eqref{new-d3} which 
corresponds to the case of dimension $d=3$. 
The other cases are similar and we leave the details to the reader.

Define 
\be{def-E1}
E_n(\varepsilon)= 
\{|\partial \kR_n|-\E[|\partial \kR_n|]\le - \varepsilon n\}, 
\ee 
and 
\be{def-E2} 
F_n(\alpha, A, \varepsilon) 
=\{\exists \Lambda \subset \Z^3\, :\, \ell_n(\Lambda)\ge \alpha n, \, 
\kI_d(\Lambda)\le A,\, \frac{n}{A\varepsilon}
\le  |\Lambda|\le \frac{A n}{\varepsilon}\}. 
\ee
We have to prove the existence of $\alpha\in (0,1)$ and $A>0$, such that for any $\varepsilon \in (0,\nu_3/2)$, 
\be{goalEF}
\pp[F_n(\alpha,A,\varepsilon)^c\mid E_n(\varepsilon)] \to 0,
\ee 
as $n\to \infty$. 
Our strategy for this is to define five sets $E_{n,1},\dots,E_{n,5}$ with $E_{n,1}=E_n(\varepsilon)$, $E_{n,5} = F_n(\alpha,A,\varepsilon)$, and 
\be{def-strategy}
\pp[ E_{n,j}\cap E^c_{n,j+1}\big]=o\big(p_n\big) \qquad \text{for }j=1,\ldots,4,
\ee
with $p_n=\pp[E_{n,1}]=\pp[E_n(\varepsilon)]$. 
This will well prove \eqref{goalEF}, since  
\begin{equation*}
E_{n,1}\cap E_{n,5}^c\ \subset \ \bigcup_{j=1}^4 E_{n,j}\cap E_{n,j+1}^c.
\end{equation*}
Set $T= n^{5/9}$. Since $T=o(n^{2/3})$, Proposition \ref{prop-slicing} and the lower bound in \eqref{shrink-d3} show that 
\be{Enxin}
\pp\left[E_{n,1}\cap E^c_{n,2}\right]
=o(p_n)\qquad \text{with}\quad
E_{n,2}=\left\{\xi_n(T)\ge \frac{\varepsilon n}{100}\right\}.
\ee
Now take $\eta$ as in Lemma \ref{lem-d3} and define 
\be{def-radius}
r:=n^{1/4}\quad \text{and}\quad
t:=\frac{\eta \varepsilon}{100} r^3, 
\ee
so that \eqref{scale-d3} is satisfied for $n$ large enough.
Let $j_0$ be the smallest integer
such that $2^{j_0/6} \ge 25\bar \kappa_3/\kappa_3$, 
with $\bar \kappa_3$ and $\kappa_3$ as in 
\eqref{shrink-d3} and \eqref{Knalpha} respectively. 
Then \eqref{main-ineqd3} and \eqref{Knalpha} show that 
\be{XinKn}
\pp\left[E_{n,2}\cap E^c_{n,3}\right]
=o(p_n)\qquad \text{with}\quad
E_{n,3}=\bigcup_{x\in \Z^d}\{ \left|\kK_n(B(x,r),t)\right| 
\ge \frac{n}{2^{3j/2}}\}.
\ee
Observe next that with our choice \reff{def-radius} 
and $L=n/2^{3j/2}$, 
$$\frac{t}{r^2} \left(1+\frac Lt\right)^{1/3} \ge \frac{t^{2/3} L^{1/3}}{r^2} = \eta^{2/3} 2^{-j_0/2}100^{-2/3}\cdot \varepsilon^{2/3} n^{1/3}.$$
Therefore a union bound and \eqref{KnGn} show that, 
with the notation of Lemma \ref{lem-key},  
\be{KnGnpn}
\pp\left[E_{n,3}\cap E^c_{n,4}\right]
=o(p_n)\qquad \text{with}\quad
E_{n,4}=\kG_n(r,t,[\delta L/t]),
\ee
where $\delta$ is the constant associated to $K= 2\cdot 100^{2/3}\bar \kappa_3\eta^{-2/3} 2^{j_0/2}$. 
Finally let $\alpha = (\delta/2)\cdot 2^{-3j_0/2}$, so that $t[\delta L/t] \ge \delta L/2= \alpha n$, for $n$ large enough. 
Note that if $\Lambda$ is the union of $m_\delta=[\delta L/t]$ disjoint 
balls of radius $r$, then at least for $n$ large enough, we have
$$
m_\delta (r/\sqrt 3)^3\ \le\  |\Lambda|\ \le \  
m_\delta (3r)^3\ \le\  \frac{10^4 \alpha}{\eta} \, \frac{n}{\varepsilon}.
$$
Therefore, if we take 
$$A=\max\left(\frac{2^{j_0+1} \sqrt 3 \, \bar \kappa_3}{ \delta^{1/3} \eta^{2/3}\kappa},\, \frac{10^4 \alpha}{\eta}\right)  ,$$ 
with $\kappa$ as in Proposition \ref{prop.cap}, then this proposition and a union bound show that $\pp(E_{n,4}\cap E_{n,5}^c)=o(p_n)$, using also that $L/t= \kO(n^{1/4})$, to remove all the combinatorial terms, and \eqref{Aminebound} for the lower bound on the volume of $\Lambda$. 
Combining this with \eqref{Enxin}, \eqref{XinKn} and \eqref{KnGnpn} 
gives \reff{def-strategy}, and concludes the proof of \reff{new-d3}.
\hfill $\square$

\section{Proof of Theorem \ref{theo-poly}}\label{sec-poly}

We start with the proof of the first statement. 
Define the positive part of a real $x$ as 
$x_+=\max(x,0)$ and $x_-=x_+-x=max(-x,0)$. Define now
\[
Z_n^{-}(\beta) = \E\left[\exp\left
(-\frac{\beta}{n^{2/d}} \, 
\big(|\partial \kR_n|-\E[|\partial \kR_n|]\big)_+\right)\right],
\]
and
\[
Z_n^{+}(\beta) = \E\left[\exp\left(\frac{\beta}{n^{2/d}} 
\, \big(|\partial \kR_n|-\E[|\partial \kR_n|]\big)_-\right)\right].
\]
Since $e^{-x}=e^{x_-}+e^{-x_+}-1$, 
we have the following simple relations (using also Jensen's
inequality)
\[
Z_n(\beta) = Z_n^+(\beta) + Z_n^-(\beta)-1,\qquad
0<Z_n^-(\beta) \le 1,\quad\text{and}\qquad
1\le Z_n(\beta)\le  Z_n^+(\beta).
\]
Thus 
$$
F^+(\beta) = \limsup\frac{1}{n^{1-2/d}} \, \log  Z_n^+(\beta)\quad \text{and}\quad F^-(\beta) = \liminf \frac{1}{n^{1-2/d}} \, \log  Z_n^+(\beta).$$
The result follows since $Z_n^+$ is 
nondecreasing in $\beta$ by construction. 

For the second statement, note first that 
$$Z_n(\beta) \ge \exp(\frac{\beta \nu_d}{4}\, n^{1-2/d})\  \pp\left[|\partial \kR_n|-\E[|\partial \kR_n|]\le - \frac{\nu_d}{4} n \right].$$
Then the fact that $\beta_d^-$ is finite follows from the lower bounds in \eqref{shrink-d3}, \eqref{shrink-d5} and \eqref{shrink-d4}. 
On the other hand, using this time the upper bounds in the latter inequalities, one get  
\begin{equation*}
\begin{split}
Z_n(\beta) \, &\le\, e^{\beta \varepsilon n^{1-2/d}}+ \sum_{k= 1}^{|\log_2 \varepsilon|} \, \exp(\beta 2^{-k+1} n^{1-2/d}) \ \pp\left[|\partial \kR_n|-\E[|\partial \kR_n|]\le - 2^{-k} n \right] \\
&\le \, e^{\beta \varepsilon n^{1-2/d}}+ \sum_{k= 1}^{|\log_2 \varepsilon|} \, \exp\Big((\beta 2^{-k+1} - \sous \kappa_d   \frac{2^{-2k/3}}{\sqrt{k\log 2}} )\cdot n^{1-2/d}\Big).
\end{split}
\end{equation*}
So we see that for $\beta$ small enough, one has for all $\varepsilon >0$, 
$$Z_n(\beta) \, \le\, e^{ \varepsilon n^{1-2/d}}+C|\log \varepsilon | \exp(-c\varepsilon^{3/4} n^{1-2/d}),$$
for some positive constants $c$ and $C$. The fact that $\beta_d^+$ is positive follows.

Let us prove the third statement now. 
So let $\beta<\beta_d^+$, $\alpha\in (0,1)$, and $A>0$ be given. Since $r_n\ge n^{2/d} (\log n)^2$, Proposition \ref{prop.cap} shows (see the discussion after its statement in the introduction, in particular \eqref{ineq-cap1}) that there exists a positive constant 
$c=c(\alpha, A)$, such that 
$$
\pp\Big[\exists \kC\in \kA(r_n,A n)\, :\, \ell_n(B(x,r_n))  \ge \alpha\frac{n}{|\kC|} \text{ for all }x\in \kC\Big]\le
\exp(-c n^{1-2/d}), 
$$ 
at least for $n$ large enough. Then, we have 
\be{beta<betac}
\Q_n^\beta\left[\exists \kC\in \kA(r_n,A n) : \ell_n(B(x,r_n))  \ge \alpha\frac{n}{|\kC|} \, \forall x\in \kC   \right]\le 
 \frac{e^{-\frac c2 n^{1-\frac 2d}}}{Z_n(\beta)} + \Q_n^\beta \left[\overline{|\partial \kR_n|}\le - \frac{c}{2\beta} n \right].
\ee
The first term on the right-hand side goes to $0$, since $Z_n(\beta)\ge 1$. For the second term, define 
$\beta'=(\beta+\beta_d^+)/2$, and note that for any $\varepsilon >0$,  
\be{polymer-meas1}
\begin{split}
\Q_n^\beta \left[\overline{|\partial \kR_n|}]\le  - \varepsilon n \right]
\le & \frac{Z_n(\beta')}{Z_n(\beta)}\cdot 
\Q_n^{\beta'}[\ind_{\{\overline{ |\partial \kR_n|}\le - \varepsilon n \}}
\exp\big(-\frac{(\beta-\beta')}{n^{2/d}}\overline{|\partial \kR_n|}\big)]\\
\le &\frac{Z_n(\beta')}{Z_n(\beta)}\cdot
\exp\left(-(\beta'-\beta) \varepsilon n^{1-2/d}\right) .
\end{split}
\ee
Then by using that $\beta'<\beta_d^+$, 
we deduce that this last term goes to $0$ as $n$ tends to infinity. 

We prove now the fourth statement. 
Fix $\beta>\beta_d^-$. By definition, there exists $c>0$, such that 
for $n$ large enough, 
\be{beta>betac1}
Z_n(\beta)\ \ge\ \exp(c n^{1-2/d}).
\ee
It follows that 
as $n\to \infty$, with $\varepsilon =c/(2\beta)$, we have
\be{beta>betac2}
\Q_n^\beta\left[|\partial \kR_n|-\E[|\partial \kR_n|]
\ge -\varepsilon n \right]\to 0.
\ee
Now, the arguments used for the proof of 
\eqref{new-d3} in the previous section   
show that there exist $\alpha$ and $A$, such that  
\begin{equation*}
\begin{split}
\pp\Big[\overline{|\partial \kR_n|} & \le -\varepsilon n,\ \min_{x\in \kC} \ell_n(B(x,r_n))< \frac{\alpha\varepsilon n}{|\kC|}\, \text{for all } \kC \in \kA(r_n,A n) \text{ satisfying }  \kI_d(B(\kC,r_n))\le A \Big] \\
& \le \exp(-2\beta n^{1-\frac 2d}), 
\end{split}
\end{equation*}
for $n$ large enough (in the proof of the previous section, take for instance $T=T(n)=(n/\sqrt{\log n})^{2/d}$, so that $r_n^d=o(T^{d/2})$). Using again that $Z_n(\beta)\ge 1$, this implies
that 
$$
\Q_n^\beta\Big[\min_{x\in \kC} \ell_n(B(x,r_n))< \frac{\alpha\varepsilon n}{|\kC|}\, \text{for all } \kC \in \kA(r_n,A n) \text{ satisfying }  \kI_d(B(\kC,r_n))\le A    \Big]\to 0, 
$$
as $n\to \infty$, and this proves \eqref{poly-3}. The proof of \eqref{poly-4} is similar and left to the reader.

It remains now to prove the last statement. 
In fact as the proof above shows, 
it suffices to see that for any fixed $\chi \in (0,1)$, \eqref{beta>betac1} holds true with 
$c=c(\beta) = (1-\chi) \beta \nu_d$, for all $\beta$ large enough (since then the proof of the previous statement works as well with $\varepsilon = (1-\chi) c/\beta$). 
For this, notice 
that \eqref{DV-inf} shows that there exists a constant $C>0$, such that 
for $n$ large enough, 
\be{polymer-meas4}
\pp\left[|\partial \kR_n|-\E[|\partial \kR_n|]
\le -(1-\chi/2)\nu_d n \right]\ \ge\  \exp(-Cn^{1-2/d}), 
\ee
as if the walk spends the first $n$ steps in a ball of volume $(\chi/4)\nu_d n$, then using also \eqref{mean-transient} we deduce that $\overline{|\partial \kR_n|}$ is smaller than $(\chi/2-1)\nu_d n$, for $n$ large enough. Therefore, for $\beta$ large enough
$$
Z_n(\beta) \ \ge \  \exp\left((\beta (1-\chi/2) \nu_d -C)n^{1-2/d}\right)\ \ge\   \exp\left(\beta (1-\chi) \nu_d n^{1-2/d}\right), 
$$
which concludes the proof.  \hfill $\square$

\vspace{0.2cm}
\noindent{\bf Acknowledgements.}
A.A. received support of the A$^*$MIDEX grant
(ANR-11-IDEX-0001-02) funded by the French Government
"Investissements d'Avenir" program.

\end{document}